\documentclass[a4paper, 12pt]{article}

\usepackage{mathrsfs}
\usepackage{epsfig}
\usepackage{amsmath}
\usepackage{amssymb}
\usepackage{latexsym}
\usepackage{amsfonts}
\usepackage{amsthm}
\usepackage[numbers,sort&compress]{natbib}
\usepackage{graphicx}
\ifx\pdfoutput\undefined  \else
 \fi
\usepackage[centerlast]{caption2}
\usepackage{color}

\newcommand {\cross } {\nu_D }

\newtheorem{definition}{Definition}

\newtheorem{lemma}{Lemma}
\newtheorem{theorem}{Theorem}
\newtheorem{observation}{Observation}
\newtheorem{claim}{Claim}
\numberwithin{figure}{section}
\numberwithin{definition}{section}
\numberwithin{observation}{section}
\numberwithin{lemma}{section}
\numberwithin{theorem}{section}
\numberwithin{table}{section}

\begin{document}

\title{
{The crossing number of cubes with small order} \footnote{The research are supported by NSFC
(11301381, 11001035, 11226280, 61303023), Science and Technology Development Fund of Tianjin Higher Institutions (20121003), Doctoral Fund of Tianjin Normal University (52XB1202).}
\author{
Guoqing Wang$^a$\thanks{Email :
gqwang1979@aliyun.com} \ \ \ \ \ \ Haoli Wang$^b$\thanks{Corresponding author's email :
bjpeuwanghaoli@163.com} \ \ \ \ \ \ Yuansheng Yang$^c$\thanks{Email: yangys@dlut.edu.cn}\\
\\
{$^a$Department of Mathematics}\\ {Tianjin Polytechnic University, Tianjin, 300387, P. R. China}\\
\\
{$^b$College of Computer and Information Engineering}\\ {Tianjin Normal University, Tianjin, 300387, P. R. China}\\
\\
{$^c$Department of Computer Science}\\ {Dalian University of Technology, Dalian, 116024, P. R. China}\\
}}

\date{}
\maketitle
\begin{abstract}

The {\it crossing number} of a graph $G$ is the minimum number of
pairwise intersections of edges in a drawing of $G$. In this paper,
we give the exact values of crossing numbers for some variations of
hypercube with order at most four, including crossed cube, locally
twisted cube and M\"{o}bius cube.
\\

\noindent {\bf Keywords:} {\it Drawing}; {\it Crossing number}; {\it
Crossed cube}; {\it Locally twisted cube}; {\it M\"{o}bius cube}
\end{abstract}

\section{Introduction}

\indent \indent  The {\it crossing number} $cr(G)$ of a graph $G$ is
the minimum number of pairwise intersections of edges in a drawing
of $G$ in the plane. The notion of crossing number is a central one
for Topological Graph Theory and has been studied extensively  by
mathematicians including Erd\H{o}s, Guy, Tur\'{a}n and Tutte, et al.
(see
\cite{ACNS82,BR80,EG73,Guy60,Kainen72,Kleitman70,WB78,Turan77,Tutte70}).
In the past thirty years, it turned out that crossing number played
an important role not only in various fields of discrete and
computational geometry (see \cite{Bi91,M02,Sz97,SoTaTo02}), but also
in VLSI theory and wiring layout problems (see
\cite{BL84,S05,L81,L83}). For this reason, the study of crossing
number of some popular parallel network topologies such as hypercube
and its variations which have good topological properties and
applications in VLSI theory, would be of theoretical importance and
practical value. Among all the popular parallel network topologies,
hypercube is the first to be studied (see
\cite{DR95,FFSV08,Madej91,SV93}). An $n$-dimensional hypercube $Q_n$
is a graph in which the nodes can be one-to-one labeled with 0-1
binary sequences of length $n$, so that the labels of any two
adjacent nodes differ in exactly one bit. It has many appealing
features such as node and edge symmetry. To improve the efficiency
of hypercube networks, several variations of hypercubes such as
crossed cube \cite{E91}, locally twisted cube \cite{YEM05} and
M\"{o}bius cube \cite{CL95} have been proposed and investigated.

Determining the crossing number of an arbitrary graph is proved to
be NP-complete \cite{GJ83}. In most cases, it is easy to find a
sufficiently ``nice" drawing for a particular kind of graph in which
the number of crossings can hardly be decreased, but is very
difficult to prove that such a drawing indeed has the smallest
possible number of crossings. Thus, it is not surprising that the
exact crossing numbers are known for graphs of few families and that
the arguments often strongly depend on their structure (see for
example
\cite{EHK81,Fiorini86,Guy72,Lin09,Yuangsheng04,Zarankiewicz54}).
With respect to cubes, the only known exact values of crossing
numbers are $cr(Q_1)=cr(Q_2)=cr(Q_3)=0$ and $cr(Q_4)=8$ \cite{DR95}.
Towards this direction,  in this paper we give the exact values of
crossing numbers for crossed cube, locally twisted cube and
M\"{o}bius cube of order at most four.

The rest of this paper is organized as follows. In Section 2 we
introduce some technical notations and tools. The crossing numbers
of crossed cube, locally twisted cube and M\"{o}bius cube of order
four are determined in Section 3, Section 4 and Section 5,
respectively.

\section{Notations and tools}



\indent \indent Let $G$ be a simple connected graph with vertex set
$V(G)$ and edge set $E(G)$. For $S\subseteq V(G)$, let $\langle
S\rangle$ be the subgraph of $G$ induced by $S$. Let $P_n$ be the
\textit{path} with $n$ vertices and let $C_n$ be the \textit{circle}
with $n$ vertices. Let $X$ and $Y$ be sets of vertices (not
necessarily disjoint) of a graph $G$. We denote by $E[X, Y]$ the set
of edges of $G$ with one end in $X$ and the other end in $Y$, and by
$e(X, Y)$ their number. If $Y = X$, we simply write $E(X)$ and
$e(X)$ for $E[X,X]$ and $e(X,X)$, respectively. When $Y =
V(G)\setminus X$, the set $E[X, Y]$ is denoted by $\partial(X)$. The
{\it degree} of a vertex $v$ in a graph $G$, denoted by $d_G(v)$, is
the number of edges of $G$ incident with $v$. When it is
unambiguous, $d_G(v)$ is abbreviated to $d(v)$.

A drawing of $G$ is said to be a {\it good} drawing, provided that
no edge crosses itself, no adjacent edges cross each other, no two
edges cross more than once, and no three edges cross in a point. It
is well known that the crossing number of a graph is attained only
in {\it good} drawings of the graph. So, we always assume that all
drawings throughout this paper are good drawings. For a good drawing
$D$ of a graph $G$, let $\nu(D)$ be the number of crossings in $D$.
In a drawing $D$, if an edge is not crossed by any other edge, we
say that it is {\it clean} in $D$.

In this paper, we will often use the term ``\textit{region}" also in
nonplanar drawings. In this case, crossings are considered to be
vertices of the ``map". The two open sets into which a simple closed
curve $C$ partitions the plane are called the {\it interior} and the
{\it exterior} of $C$. By a \textit{line segment}, we mean a curve
incident with vertices or crossings. The \textit{boundary} of a
region $f$ is the boundary of the open set $f$ in the usual
topological sense. A region is said to be \textit{incident} with the
vertices and line segaments in its boundary.

Two regions are \textit{adjacent} if their boundaries have a line
segment in common. Let $f,h$ are two regions of a graph $G$. Let
$$\mathcal {B}(f,h)=\left\{
\begin{array}{llll}
1, &\hbox{if $f,h$ are adjacent,}\\
0, &\hbox{otherwise.}
\end{array}
\right.$$ Let $G_1$ and  $G_2$ be two connected graphs, and let
$D_1$ be a good drawing of $G_1$. For a region $f$ of $G_1$ in the
drawing $D_1$, we define
$$V_{in}(f;G_2)=\{v\in V(G_2): v \mbox{ lies in the region } f\}$$
and $$V_{on}(f)=\{v\in V(G_1):v \mbox{ is incident with the region }
f\}.$$ Two drawings of $G$ are \textit{isomorphic} if and only if
there is an incidence preserving one-to-one correspondence between
their vertices, edges, parts of edges and regions.

Now we give the definitions of variations of hypercubes which are
studied in this paper. For more details, one can refer to
\cite{X01}.

\begin{definition}(Crossed cube) Two binary strings $x=x_2x_1$ and $y=y_2y_1$ are {\it pair-related},
denoted by $x\sim y$, if and only if
$(x,y)\in\{(00,00),(10,10),(01,11)$, $(11,01)\}$. The
$n$-dimensional crossed cube $CQ_n$  has vertex set $V=\{x_n\cdots
x_2x_1:x_i\in \{0,1\},i=1,2,\ldots,n\}$, and two vertices
$x=x_n\cdots x_2x_1$ and $y=y_n\cdots y_2y_1$ are linked by an edge
if and only if there exists $j \ \ (1 \leq j \leq n)$ such
that \\
\indent(a) $x_n\cdots x_{j+1}=y_n\cdots y_{j+1}$, \\
\indent(b) $x_j\neq y_j$, \\
\indent(c) $x_{j-1}=y_{j-1}$ if $j$ is even, and \\
\indent(d) $x_{2i}x_{2i-1}\sim y_{2i}y_{2i-1}$ for each
$i=1,2,\ldots,\lceil \frac{j}{2}\rceil-1$.
\end{definition}

The graphs shown in Figure 2.1 are $CQ_1$, $CQ_2$, $CQ_3$ and
$CQ_4$, respectively.
\begin{figure}[ht]
\centering
\includegraphics[scale=1.1]{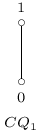} \hspace{20pt}
\includegraphics[scale=1.1]{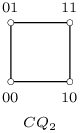} \hspace{20pt}
\includegraphics[scale=1.1]{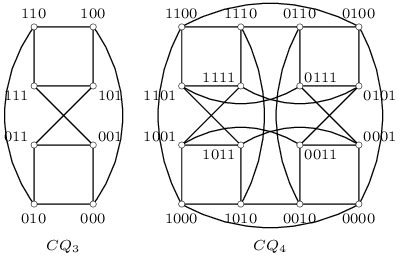}
\caption{\small{Crossed cubes $CQ_1$, $CQ_2$, $CQ_3$ and $CQ_4$}}
\end{figure}

\begin{definition}(Locally twisted cube)
The {\it $n$-dimensional locally twisted cube} $LTQ_n(n\geq 2)$ is defined recursively as follows. \\
\indent (a) $LTQ_2$ is a graph isomorphic to $Q_2$. \\
\indent (b) For $n\geq 3$, $LTQ_n$ is built from two disjoint copies
of $LTQ_{n-1}$ according to the following steps. Let $0LTQ_{n-1}$
denote the graph obtained by prefixing the label of each vertex of
one copy of $LTQ_{n-1}$ with 0, let $1LTQ_{n-1}$ denote the graph
obtained by prefixing the label of each vertex of the other copy
$LTQ_{n-1}$ with 1, and connect each vertex $x=0x_2x_3\ldots x_n$ of
$0LTQ_{n-1}$ with the vertex $1(x_2+x_n)x_3\ldots x_n$ of
$1LTQ_{n-1}$ by an edge, where $‘+’$ represents the modulo 2
addition.
\end{definition}

The graphs shown in Figure 2.2 are $LTQ_2$, $LTQ_3$ and $LTQ_4$,
respectively.
\begin{figure}[ht]
\centering
\includegraphics[scale=1.1]{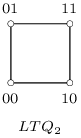} \hspace{20pt}
\includegraphics[scale=1.1]{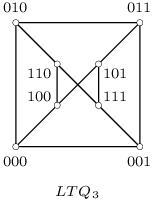} \hspace{20pt}
\includegraphics[scale=1.1]{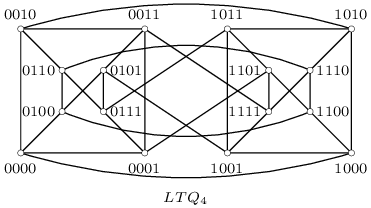}
\caption{\small{Locally twisted cubes $LTQ_2$, $LTQ_3$ and $LTQ_4$}}
\end{figure}

\begin{definition}(M\"{o}bius cube)
The $n$-dimensional M\"{o}bius cube $MQ_n$ is such a graph with
vertex set $V=\{x_1x_2\cdots x_n:x_i\in \{0,1\},i=1,2,\ldots,n\}$,
and the vertex $X=x_1x_2\cdots x_n$ connects to $n$ other vertices
$Y_i(1\leq i\leq n)$, where each $Y_i$ satisfies one of the following equations:
$$Y_i=\left\{
\begin{array}{llll}
x_1\ldots x_{i-1}\bar{x}_ix_{i+1}\ldots x_n, &\hbox{if $x_{i-1}=0$,}\\
x_1\ldots x_{i-1}\bar{x}_i\bar{x}_{i+1}\ldots \bar{x}_n, &\hbox{if
$x_{i-1}=1$,}
\end{array}
\right.$$ where $\bar{x}_i$ is the complement of $x_i$ in \{0,1\}.
\end{definition}
From the above definition, $X$ connects to $Y_i$ by complementing
the bit $x_i$ if $x_{i-1}=0$ or by complementing all bits of
$x_i,\cdots,x_n$ if $x_{i-1}=1$. The connection between $X$ and
$Y_1$ is undefined, so we can assume $x_0$ is either equal to 0 or
equal to 1, which gives us slightly different network topologies. If
we assume $x_0=0$, we call the network a ``0-M\"{o}bius cube"; and
if we assume $x_0=1$, we call the network a ``1-M\"{o}bius cube",
denoted by $0\mbox{-}MQ_n$ and $1\mbox{-}MQ_n$, respectively. The
graphs shown in Figure 2.3 and Figure 2.4 are 0-$MQ_i$ and 1-$MQ_i$
for $i=1,2,3,4$.
\begin{figure}[ht]
\centering
\includegraphics[scale=1.0]{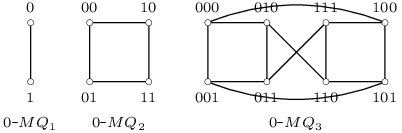} \hspace{10pt}
\includegraphics[scale=1.0]{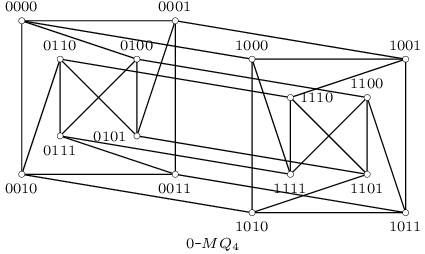}
\caption{\small{M\"{o}bius cubes 0-$MQ_1$, 0-$MQ_2$, 0-$MQ_3$ and
0-$MQ_4$}}
\end{figure}
\begin{figure}[ht]
\centering
\includegraphics[scale=1.0]{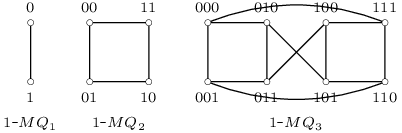} \hspace{10pt}
\includegraphics[scale=1.0]{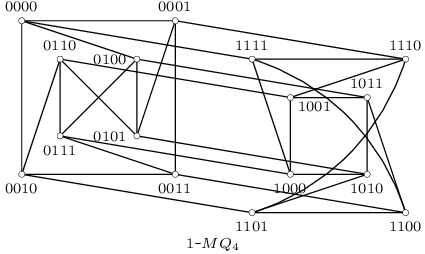}
\caption{\small{M\"{o}bius cubes 1-$MQ_1$, 1-$MQ_2$, 1-$MQ_3$ and
1-$MQ_4$}}
\end{figure}

We next present some tools which will be useful for later sections
of this paper.

Let $A$ and $B$ be two disjoint subsets of $E$. In a drawing $D$,
the number of the crossings formed by an edge in $A$ and another
edge in $B$ is denoted by $\cross(A,B)$. The number of the crossings
that involve a pair of edges in $A$ is denoted by $\cross (A)$. Then
$\nu(D) = \cross(E)$. Then the following lemma is straightforward.

\begin{lemma} \label{lemma counting the number of crossings}
Let $A$, $B$, $C$ be mutually disjoint subsets of $E$. Then,
$$\begin {array}{llll}
\cross(C,A\cup B)&=&\cross(C,A)+\cross(C,B),\\
\cross(A\cup B)&=&\cross(A)+\cross(B)+\cross(A,B).
\end {array}$$
\end{lemma}

In 1995, A.M. Dean and R.B. Richter proved the following.

\begin{lemma} \label{lemma a good drawing about a vertex with degree 4}\cite{DR95}
Let $v$ be a vertex with degree 4 in a graph $G$. In an optimal
drawing of $G$, there is no edge that crosses exactly three of the
edges incident with $v$.
\end{lemma}

It is easy to see that $CQ_3\cong LTQ_3\cong 0\mbox{-}MQ_3\cong
1\mbox{-}MQ_3$.

\bigskip

 \textsf{Throughout this section, let $G_3$
be any one of the above isomorphic graphs $CQ_3$, $LTQ_3$,
$0\mbox{-}MQ_3$ and $1\mbox{-}MQ_3$.}

\begin{observation}\label{Observation G3 contains K3,3}
$G_3$ contains a subgraph homeomorphic to $K_{3,3}$.
\end{observation}

\begin{observation}\label{Observation two vertex-disjoint cycles}
Let $C$ and $C^{'}$ be two vertex-disjoint cycles. Then the number
of crossings between $E(C)$ and $E(C^{'})$ is even.
\end{observation}

\begin{lemma}\label{Lemma |Von(f)| leq 6}
Assume that $V(G_3)=V_1\cup V_2$  such that $\langle
V_1\rangle\cong\langle V_2\rangle\cong C_4$. If all the edges of
$\langle V_1\rangle$ are clean, then $|V_{on}(f)|\leq 6$ for every
region $f$ of $G_3$.
\end{lemma}

\begin{proof} Let $\langle V_1\rangle=u_1u_2u_3u_4u_1$. Since all the edges of $\langle
V_1\rangle$ are clean, we may assume $\langle V_2\rangle $ lie in
the outer space of cycle $u_1u_2u_3u_4u_1$. Let $h$ be the inner
space of cycle $u_1u_2u_3u_4u_1$. Let $f_i$ be the region of $G_3$,
distinct from $h$, such that the edge $u_iu_{i+1}$ is on the
boundary of $f_{i}$ for $i=1,2,3,4$, where the subscripts are taken
modulo $4$. It is easy to see that
\begin{equation}\label{equation the region fi}
u_{i+2},u_{i+3}\not\in V_{on}(f_i)
\end{equation} for $i=1,2,3,4$.

Suppose to the contrary that there exists a region $f$ of $G_3$ such
that $|V_{on}(f)|\geq 7$. Then at least three of
$\{u_1,u_2,u_3,u_4\}$, say $u_1,u_2,u_3\in V_{on}(f)$. It follows
that at least two edges of $u_1u_2,u_2u_3,u_3u_4,u_4u_1$ are on the
boundary of $f$. Since $|V_{on}(f)|=7$, we have that $f$ is $f_i$
for some $i\in \{1,2,3,4\}$. Then we derive a contradiction with
\eqref{equation the region fi}.
\end{proof}

Before proceeding, we need the following definition.


\begin{lemma} \label{Lemma the number of emitting edges}
Let $X$ be a subset of $V(G_3)$ with $1\leq |X|\leq 4$, and let
$X'=V(G_3)\setminus X$. Then
$$\begin{array}{llll} |\partial(X)|\geq \left\{\begin{array}{llll}
               3,  & \mbox{if \ \ } |X|=1;\\
               6,  & \mbox{if \ \ } |X|=2 \mbox{ and } \langle X\rangle \not\cong P_2;\\
               5,  & \mbox{if \ \ } |X|=3;\\
               6,  & \mbox{if \ \ } |X|=4 \mbox{ and } \langle X\rangle \not\cong C_4;\\
               4,  & \mbox{Otherwise.} \\
              \end{array}
              \right .
\end{array}$$
Moreover, if $|X|=3$ and $|\partial(X)|=5$, then $\langle
X'\rangle\cong C_5$ or $\langle X'\setminus \{v\}\rangle\cong C_4$
for some $v\in X'$.
\end{lemma}

\begin{proof} Since $G_3$ is a 3-regular graph, we have
\begin{equation}\label{equation delta X}
|\partial(X)|=\sum\limits_{v\in X}d(v)-2e(X)=3|X|-2e(X).
\end{equation}

Suppose $|X|=1$. By \eqref{equation delta X}, it is clear that
$|\partial(X)|=3$.

Suppose $|X|=2$. Since $e(X)\leq 1$, it follows from \eqref{equation
delta X} that $|\partial(X)|\geq 3\times 2-2=4$. Moreover, if
$\langle X\rangle \not\cong P_2$, then $e(X)=0$. It follows from
\eqref{equation delta X} that $|\partial(X)|\geq 3\times 2=6$.

Suppose $|X|=3$.  Since $G_3$ is triangle-free, we have $e(X)\leq
2$. It follows from \eqref{equation delta X} that $|\partial(X)|\geq
3\times 3-2\times 2=5$. If $|\partial(X)|=5$, we have $|\partial
(X)|=|\partial(X')|=3|X'|-2e(X')=3\times 5-2e(X')$, which implies $e
(X')=5$. Since $G_3$ is triangle-free, it is not hard to infer that
$\langle X'\rangle\cong C_5$ or $\langle X'\setminus
\{v\}\rangle\cong C_4$ for some $v\in X'$.

Suppose $|X|=4$. Since $G_3$ is triangle-free, we have $e(X)\leq 4$.
It follows from \eqref{equation delta X} that $|\partial(X)|\geq
3\times 4-2\times 4=4$. Moreover, if $\langle X\rangle \not\cong
C_4$, then $e(X)\leq 3$. It follows from \eqref{equation delta X}
that $|\partial(X)|\geq 3\times 4-2\times 3=6$.
\end{proof}

\begin{lemma}\label{Lemma three vertex set (6 1 1)}
Let $X_1,X_2,X_3$ be pairwise disjoint vertex subset of $G_3$ with
$V(G_3)=X_1\cup X_2\cup X_3$ and $|X_1|\geq |X_2|\geq |X_3|>0$. Then
$\sum\limits_{1\leq i<j\leq 3}e(X_i,X_j)\geq 5$, and moreover, the
equality implies $(|X_1|,|X_2|,|X_3|)=(6,1,1)$ and $e(X_i,X_j)>0$
for any $1\leq i<j\leq 3$.
\end{lemma}

\begin{proof} Since $|X_1|\geq 2$, it follows from Lemma \ref{Lemma the number
of emitting edges} that $$\sum\limits_{1\leq i<j\leq
3}e(X_i,X_j)=\frac{1}{2}\sum\limits_{1\leq i\leq
3}\sum\limits_{\stackrel{1\leq j\leq 3}{j\neq i}}e(X_i,X_j)\geq
\frac{1}{2}\times (4+3+3)=5.$$ Suppose $\sum\limits_{1\leq i<j\leq
3}e(X_i,X_j)=5$. Then $\sum\limits_{\stackrel{1\leq j\leq 3}{j\neq
i}}e(X_i,X_j)=3$ for $i=2,3$, and so $|X_2|=|X_3|=1$, i.e.,
$(|X_1|,|X_2|,|X_3|)=(6,1,1)$. Moreover, since $\sum\limits_{1\leq
i<j\leq 3}e(X_i,X_j)\geq |\partial(X_s)|+
|\partial(X_t)|-e(X_s,X_t)\geq 3+3-e(X_s,X_t)$, we have
$e(X_s,X_t)>0$, where $1\leq s<t\leq 3$.
\end{proof}

\begin{lemma}\label{lemma three vertex set equal 6}
Let $X_1,X_2,X_3$ be pairwise disjoint vertex subset of $G_3$ with
$V(G_3)=X_1\cup X_2\cup X_3$ and $|X_1|\geq |X_2|\geq |X_3|>0$.
Suppose that $\sum\limits_{1\leq i<j\leq 3}e(X_i,X_j)=6$, and that
$e(X_i,X_j)>0$ for all $1\leq i<j\leq 3$, and that there exist
$1\leq s< t\leq 3$ with $e(X_s,X_t)=1$.
Then one of the following conditions holds. \\
(i)  $(|X_1|,|X_2|,|X_3|)=(5,2,1)$ with $e(X_2,X_3)=1$ and $\langle
X_2\rangle\cong P_2$; \\
(ii) $(|X_1|,|X_2|,|X_3|)=(4,3,1)$ with $e(X_1,X_3)=1$ and $\langle
X_1\rangle\cong C_4$.
\end{lemma}

\begin{proof} We note first that
\begin{equation}\label{equation 6=..=Xl+Xm-e(Xl,Xm)}
3(|X_{\ell}|+|X_m|)-2(e(X_{\ell})+e(X_m))-e(X_{\ell},X_m)=6
\end{equation}
for any $1\leq \ell<m\leq 3$. This is because that $G_3$ is a
3-regular graph and that $\sum\limits_{1\leq i<j\leq
3}e(X_i,X_j)=|\partial(X_{\ell})|+|\partial(X_m)|-e(X_{\ell},X_m)$.

\textbf{Claim 1.} If $|X_3|=1$ then $t=3$.

Assume to the contrary that $|X_3|=1$ and $t\neq 3$, i.e.,
$(s,t)=(1,2)$. Since $G_3$ is a 3-regular graph, we have
$|\partial(X_3)|=3$. Since $e(X_1,X_2)=1$, it follows that
$\sum\limits_{1\leq i<j\leq
3}e(X_i,X_j)=|\partial(X_3)|+e(X_1,X_2)=3+1<6$, a contradiction.
This proves Claim 1. \qed

We observe that all the possible cases of $(|X_1|,|X_2|,|X_3|)$ are
$(6,1,1),(5,2,1),(4,3,1),$ $(4,2,2),(3,3,2)$.

Suppose $(|X_1|,|X_2|,|X_3|)=(6,1,1)$. Applying \eqref{equation
6=..=Xl+Xm-e(Xl,Xm)} with $\ell=2$ and $m=3$, since
$e(X_{\ell},X_m)>0$, we derive a contradiction. Hence,
$$(|X_1|,|X_2|,|X_3|)\neq (6,1,1).$$

Applying \eqref{equation 6=..=Xl+Xm-e(Xl,Xm)} with $\ell=s$ and
$m=t$, since $e(X_s,X_t)=1$, we derive that
\begin{equation}\label{equation |Xs|+|Xt| is even}
|X_s|+|X_t|\equiv 1\pmod 2.
\end{equation}

By \eqref{equation |Xs|+|Xt| is even}, we have
$$(|X_1|,|X_2|,|X_3|)\neq (4,2,2).$$

Suppose $(|X_1|,|X_2|,|X_3|)=(3,3,2)$. By \eqref{equation |Xs|+|Xt|
is even}, we have $(|X_s|,|X_t|)=(3,2)$. Since $G_3$ is triangle
free, we have that $e(X_s)\leq 2$ and $e(X_t)\leq 1$. Since $e(X_s,
X_t)=1$, we have that the left side of \eqref{equation
6=..=Xl+Xm-e(Xl,Xm)} is no less than $3\times (3+2)-2\times
(2+1)-1=8$, which is a contradiction. Hence,
$$(|X_1|,|X_2|,|X_3|)\neq (3,3,2).$$

Suppose $(|X_1|,|X_2|,|X_3|)=(5,2,1)$. By Claim 1 and
\eqref{equation |Xs|+|Xt| is even}, we have $(s,t)=(2,3)$, i.e.,
$|X_s|=2$ and $|X_t|=1$. Applying  \eqref{equation
6=..=Xl+Xm-e(Xl,Xm)} with $\ell=s$ and $m=t$, since $e(X_s,X_t)=1$,
we have that $e(X_s)+e(X_t)=1$. Since $|X_t|=1$, we have $e(X_s)=1$.
This implies $\langle X_s\rangle\cong P_2$, done.

Suppose $(|X_1|,|X_2|,|X_3|)=(4,3,1)$. By Claim 1 and
\eqref{equation |Xs|+|Xt| is even}, we have $(s,t)=(1,3)$, i.e.,
$|X_s|=4$ and $|X_t|=1$. Applying  \eqref{equation
6=..=Xl+Xm-e(Xl,Xm)} with $\ell=s$ and $m=t$, since $e(X_s,X_t)=1$,
we have that $e(X_s)+e(X_t)=4$. Since $|X_t|=1$, we have $e(X_s)=4$.
Since $G_3$ is triangle-free, we conclude that $\langle
X_s\rangle\cong C_4$, we are done.
\end{proof}

\begin{lemma}\label{lemma three vertex set geq 8}
Let $X_1,X_2,X_3$ be pairwise disjoint vertex subset of $G_3$ with
$V(G_3)=X_1\cup X_2\cup X_3$ and $|X_1|\geq |X_2|\geq |X_3|>0$.
Suppose that there exist $1\leq s<t\leq 3$ such that $e(X_s,X_t)=0$.
Then $|X_t|=1$ or $\sum\limits_{1\leq i<j\leq 3}e(X_i,X_j)\geq 8$.
\end{lemma}

\begin{proof} Suppose $|X_t|\geq 2$, and so $|X_s|\geq |X_t|\geq 2$.
Since $e(X_s,X_t)=0$, by Lemma \ref{Lemma the number of emitting
edges}, we have $\sum\limits_{1\leq i<j\leq
3}e(X_i,X_j)=e(X_s,X_j)+e(X_t,X_j)\geq 4+4=8$ where
$j\in\{1,2,3\}\setminus \{s,t\}$. We are done.
\end{proof}

\begin{lemma}\label{Lemma four region geq 7}
Let $X_1,X_2,X_3,X_4$ be pairwise disjoint vertex subsets of $G_3$
with $V(G_3)=X_1\cup X_2\cup X_3\cup X_4$ and $|X_1|\geq |X_2|\geq
|X_3|\geq |X_4|>0$. Then $\sum\limits_{1\leq i<j\leq
4}e(X_i,X_j)\geq 7$, and moreover, the equality
implies that one of the following conditions holds. \\
 (i) $e(X_i,X_j)>0$ for any $1\leq i<j\leq 4$;\\
 (ii) $|X_1|=5$ and there exist $2\leq \alpha<\beta\leq 4$ such that
$e(X_{\alpha},X_{\beta})=0$ and $e(X_i,X_j)>0$ for all $1\leq
i<j\leq 4$ with $(i,j)\neq (\alpha,\beta)$.
\end{lemma}

\begin{proof} Since $|X_1|\geq 2$, it follows from Lemma \ref{Lemma the number of emitting edges} that $$\sum\limits_{1\leq
i<j\leq 4}e(X_i,X_j)=\frac{1}{2}\sum\limits_{1\leq i\leq
4}\sum\limits_{\stackrel{1\leq j\leq 4}{j\neq i}}e(X_i,X_j)\geq
\frac{1}{2}\times (4+3+3+3)=6.5,$$ and so $\sum\limits_{1\leq
i<j\leq 4}e(X_i,X_j)\geq 7$.

Now assume that $\sum\limits_{1\leq i<j\leq 4}e(X_i,X_j)= 7$, and
that (i) does not hold, i.e., there exist $1\leq s<t\leq 4$ such
that
$$e(X_s,X_t)=0.$$ It suffices to show (ii) holds. Let $k,\ell\in
\{1,2,3,4\}\setminus \{s,t\}$ with $k<\ell$. Since $e(X_s,X_t)=0$,
it follows that
\begin{equation}\label{equation emiting edges of Xs and Xt}
e(X_k,X_{\ell})+\sum\limits_{\stackrel{1\leq j\leq 4}{j\neq
s}}e(X_s,X_j)+\sum\limits_{\stackrel{1\leq j\leq 4}{j\neq
t}}e(X_t,X_j)= \sum\limits_{1\leq i<j\leq 4}e(X_i,X_j)=7.
\end{equation}
By \eqref{equation emiting edges of Xs and Xt} and Lemma \ref{Lemma
the number of emitting edges}, we conclude that $$|X_t|=1,$$ and
that either $\langle X_s\rangle \cong C_4$ or $\langle X_s\rangle
\cong P_2$ or $|X_s|=1$.

Suppose $\langle X_s\rangle \cong C_4$. By Observation
\ref{observation 4-cycle match} (ii), we have $e(X_s,X_t)\neq 0$, a
contradiction.

Suppose $\langle X_s\rangle \cong P_2$. By \eqref{equation emiting
edges of Xs and Xt} and Lemma \ref{Lemma the number of emitting
edges}, we have $$e(X_k,X_{\ell})=0.$$ By Observation
\ref{observation 4-cycle match} (i) and (ii), we conclude that
$\langle X_s\cup X_{\ell}\rangle\cong \langle X_t\cup
X_k\rangle\cong C_4$, which implies $|X_{\ell}|=2>|X_t|$, and thus,
$e(X_k,X_{\ell})\geq 1$, a contradiction.

Suppose $|X_s|=1$. Then $(|X_k|,|X_{\ell}|)\in
\{(3,3),(4,2),(5,1)\}$. By \eqref{equation emiting edges of Xs and
Xt} and Lemma \ref{Lemma the number of emitting edges}, we have
$$e(X_k,X_{\ell})\leq 1.$$ Combined with Observation \ref{observation 4-cycle
match} (ii), we conclude $\langle X_k\rangle \not\cong C_4$. If
$(|X_k|,|X_{\ell}|)\in \{(3,3),(4,2)\}$, by Lemma \ref{Lemma the
number of emitting edges}, we have $\sum\limits_{\stackrel{1\leq
j\leq 4}{j\neq k}}e(X_k,X_j)+\sum\limits_{\stackrel{1\leq j\leq
4}{j\neq \ell}}e(X_{\ell},X_j)\geq 9$, and thus,
$$\sum\limits_{1\leq i<j\leq 4}e(X_i,X_j)\geq
\sum\limits_{\stackrel{1\leq j\leq 4}{j\neq
k}}e(X_k,X_j)+\sum\limits_{\stackrel{1\leq j\leq 4}{j\neq
\ell}}e(X_{\ell},X_j)-e(X_k,X_{\ell})\geq 8,$$ a contradiction.
Hence
$$(|X_k|,|X_{\ell}|)=(5,1).$$ Since $e(X_s,X_t)=0$ and
$7=\sum\limits_{\stackrel{1\leq j\leq 4}{j\neq k}}e(X_k,X_j) \geq
\sum\limits_{\stackrel{1\leq j\leq 4}{j\neq
s}}e(X_s,X_j)+\sum\limits_{\stackrel{1\leq j\leq 4}{j\neq
\ell}}e(X_{\ell},X_j)+\sum\limits_{\stackrel{1\leq j\leq 4}{j\neq
t}}e(X_t,X_j)-(e(X_s,X_{\ell})+e(X_s,X_t)+e(X_t,X_{\ell}))$, it
follows from Lemma \ref{Lemma the number of emitting edges} that
$e(X_s,X_{\ell})=e(X_t,X_{\ell})=1$. By Lemma \ref{Lemma the number
of emitting edges}, $|\partial(X_1)|\geq 5$, which implies that
$e(X_1,X_j)>0$ for $2\leq j\leq 4$. Put $\alpha,\beta$ to be $s,t$,
then the lemma follows. \end{proof}

\begin{lemma}\label{Lemma five regions geq 8}
For $t\geq 5$, let $X_1,\ldots,X_t$ be pairwise disjoint vertex
subsets of $G_3$ with $V(G_3)=\bigcup\limits_{i=1}^t X_i$, where
$|X_1|\geq |X_2|\geq \cdots\geq |X_t|>0$. Then $\sum\limits_{1\leq
i<j\leq t}e(X_i,X_j)\geq 8$.
\end{lemma}

\begin{proof} By Lemma \ref{Lemma the number of emitting edges}, we have $\sum\limits_{1\leq i<j\leq
t}e(X_i,X_j)=\frac{1}{2}\sum\limits_{1\leq i\leq
t}\sum\limits_{\stackrel{1\leq j\leq t}{j\neq i}}e(X_i,X_j)\geq
\frac{1}{2}\times 3\times 5=7.5$, and so $\sum\limits_{1\leq i<j\leq
t}e(X_i,X_j)\geq 8$.
\end{proof}

\section{Crossing number of $CQ_4$}

\indent \indent In Figure 3.1, we give a drawing of $CQ_4$ with 8
crossings. Hence, we have the following

\begin{lemma} \label{lemma upper bound of CQ4}
$cr(CQ_4)\leq 8$.
\end{lemma}
\begin{figure}[ht]
\centering
\includegraphics[scale=0.9]{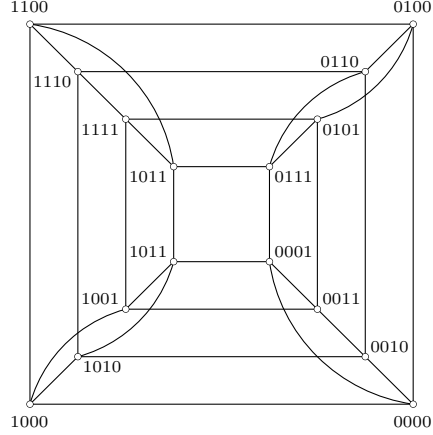}
\caption{\small{A good drawing of $CQ_4$ with 8 crossings}}
\end{figure}

In the rest of this section, we shall prove that the  value of
$cr(CQ_4)$ is exactly equal to 8.  We rename the vertices of $CQ_4$
as shown in Figure 2.8. Let $\ell\mbox{-}CQ_3$ and $r\mbox{-}CQ_3$
be the subgraphs induced by vertex subset $\{v_i:0\leq i\leq 7\}$
and by vertex subset $\{v_i:8\leq i\leq 15\}$, respectively. Note
that both $\ell\mbox{-}CQ_3$ and $r\mbox{-}CQ_3$ are isomorphic to
$CQ_3$. For convenience, we abbreviate
$$\begin {array}{llll}
E_{\ell}&=&E(\ell\mbox{-}CQ_3),\\
E_r&=&E(r\mbox{-}CQ_3),\\
E_{\ell,r}&=&E[V(\ell\mbox{-}CQ_3),V(r\mbox{-}CQ_3)].\\
\end {array}$$
\begin{figure}[ht]
\centering
\includegraphics[scale=0.9]{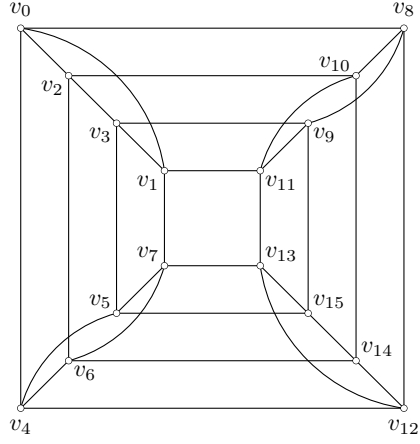}
\caption{\small{ $CQ_4$ with the renamed vertices}}
\end{figure}

By Lemma \ref{lemma counting the number of crossings}, we have the
following

\begin{lemma}\label{Lemma total crossing and left }
Let $D$ be a good drawing of $CQ_4$. Then
$$\nu(D)=\cross(E_{\ell})+\cross(E_r)+\cross(E_{\ell,r})+\cross(E_{\ell},E_r)+\cross(E_{\ell},E_{\ell,r})+\cross(E_r,E_{\ell,r}).$$
\end{lemma}

\vspace{3pt}

\begin{observation} \label{observation vertex partition of LTQ3}
Let $G$ be a graph isomorphic to $CQ_3$. Then there exist exactly
two partitions of $V(G)$: $V(G)=V_1\cup V_2$ and $V(G)=V_3\cup V_4$
such that $\langle V_i\rangle\cong C_4$ for $i=1,2,3,4$. In
particular, $\langle V_i\cap V_j\rangle=P_2$ for any $i\in \{1,2\}$
and $j\in \{3,4\}$ (see Figure 3.3(2)).
\end{observation}
\begin{figure}[ht]
\centering
\includegraphics[scale=0.9]{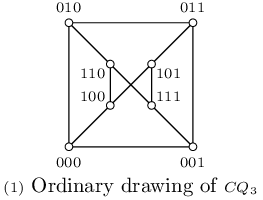} \hspace{30pt}
\includegraphics[scale=0.9]{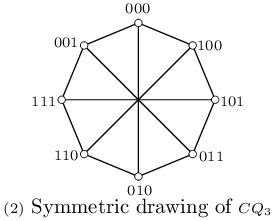}
\caption{\small{Two drawings of $CQ_3$}}
\end{figure}

\begin{theorem}\label{theorem cr(CQ4)=8}
$cr(CQ_4)=8$.
\end{theorem}
\begin{proof} By Lemma \ref{lemma upper bound of CQ4}, it suffices to prove that
$cr(CQ_4)\geq 8$. Suppose to the contrary that
\begin{equation}\label{equation cr(CQ4) leq 7}
cr(CQ_4)\leq 7,
\end{equation}
i.e., there exists a good drawing $D$ of $CQ_4$ such that
$\nu(D)\leq 7$. Without loss of generality, we may assume
$\cross(E_{\ell})+\cross(E_{\ell},E_{\ell,r})\leq
\cross(E_r)+\cross(E_r,E_{\ell,r})$.
Since $\nu(D)\leq 7$, by Lemma
\ref{Lemma total crossing and left }, we have
\begin{equation}\label{left leq 3.5(2)}
\cross(E_{\ell})+\cross(E_{\ell},E_{\ell,r})+0.5\cross(E_{\ell,r})+0.5\cross(E_{\ell},E_r)\leq
3.5.
\end{equation}
By Observation \ref{Observation G3 contains K3,3}, we have
$\cross(E_{\ell})\geq 1$. It follows from \eqref{left leq 3.5(2)}
that
\begin{equation}\label{equation 1 leq cr(El) leq 3(2)}
1\leq \cross(E_{\ell})\leq 3.
\end{equation}

By Observation \ref{observation vertex partition of LTQ3}, we may
assume without loss of generality that there exists a permutation
$a,b,c,d,a^{'},b^{'},c^{'},d^{'}$ of $V(\ell\mbox{-}CQ_3)$
satisfying that
\begin{equation}\label{equation CQ four 4-cycles}
\langle\{a,b,c,d\}\rangle\cong
\langle\{a^{'},b^{'},c^{'},d^{'}\}\rangle\cong \langle
\{a,b,a^{'},b^{'}\}\rangle\cong \langle
\{c,d,c^{'},d^{'}\}\rangle\cong C_4,
\end{equation}
 and moreover,
\begin{equation}\label{equation CQ four 4-cycles in particular 1}
\langle\{a,b,c,d\}\rangle=abcda,
\end{equation}
\begin{equation}\label{equation CQ four 4-cycles in particular 2}
\langle\{a^{'},b^{'},c^{'},d^{'}\}\rangle=a^{'}b^{'}c^{'}d^{'}a^{'}.
\end{equation}

\noindent\textbf{Claim 1.} {\sl If $\cross(E_{\ell})=3$ then there
exists a region $f$ of $\ell\mbox{-}CQ_3$ with $|V_{on}(f)|=8$, and
if $\cross(E_{\ell})=2$ then there exists a region $f$ of
$\ell\mbox{-}CQ_3$ with $|V_{on}(f)|\geq 7$.}

{\sl Proof of Claim 1.} Suppose $\cross(E_{\ell})=3$. By \eqref{left
leq 3.5(2)} and Lemma \ref{Lemma the number of emitting edges}, we
conclude that $\cross(E_{\ell},E_{\ell,r})=0$ and
$\cross(E_{\ell},E_r)=0$, i.e., all vertices of $r\mbox{-}CQ_3$ lie
in the same region of $\ell\mbox{-}CQ_3$. Since
$\cross(E_{\ell},E_{\ell,r})=0$, it follows that there exists a
region $f$ of $\ell\mbox{-}CQ_3$ with $|V_{on}(f)|=8$.

Suppose that $\cross(E_{\ell})=2$ and
\begin{equation}\label{equation V_on(f) leq 6(2)}
|V_{on}(f)|\leq 6
\end{equation}
for every region $f$ of $\ell\mbox{-}CQ_3$. By \eqref{left leq
3.5(2)}, we have that
\begin{equation}\label{equation cr(el,elr) leq 1(2)}
\cross(E_{\ell},E_{\ell,r})\leq 1
\end{equation}
and
\begin{equation}\label{equation cr(el,er) leq 3(2)}
\cross(E_{\ell},E_r)\leq 3
\end{equation}
If $\cross(E_{\ell},E_r)=0$, i.e., all vertices of $r\mbox{-}CQ_3$
lie in the same region of $\ell\mbox{-}CQ_3$, by \eqref{equation
V_on(f) leq 6(2)}, we have $\cross(E_{\ell},E_{\ell,r})\geq 2$, a
contradiction with \eqref{equation cr(el,elr) leq 1(2)}. By Lemma
\ref{Lemma the number of emitting edges}, we conclude that
$\cross(E_{\ell},E_r)=3$ and there exists a region $f$ of
$\ell\mbox{-}CQ_3$ such that $|V_{in}(f,\ell\mbox{-}CQ_3)|=7$, which
is a contradiction with \eqref{equation V_on(f) leq 6(2)}. This
proves Claim 1. \qed

\noindent\textbf{Claim 2.} {\sl
$\cross(E(\{a,b,c,d\}),E(\{a^{'},b^{'},c^{'},d^{'}\}))=\cross(E(\{a,b,a^{'},b^{'}\}),E(\{c,d,
c^{'},d^{'}\}))=0$.}

{\sl Proof of Claim 2.} Without loss of generality, we may assume to
the contrary that
\begin{equation}\label{equation CQ Claim 2 assume contrary}
\cross(E(\{a,b,c,d\}),E(\{a^{'},b^{'},c^{'},d^{'}\}))\neq 0.
\end{equation}
This implies that
\begin{equation}\label{equation Von(f) leq 7(2)}
|V_{on}(f)|\leq 7,
\end{equation}
where $f$ runs over all regions of $\ell\mbox{-}CQ_3$. By
\eqref{equation 1 leq cr(El) leq 3(2)} and \eqref{equation CQ Claim
2 assume contrary}, we have $$3\geq \cross(E_{\ell})\geq
\cross(E(\{a,b,c,d\}),E(\{a^{'},b^{'},c^{'},d^{'}\}))\geq 1.$$ By
Observation \ref{Observation two vertex-disjoint cycles} and Claim
1, we conclude that
\begin{equation}\label{equation cross(E0,E1)=2(2)}
\cross(E_{\ell})=\cross(E(\{a,b,c,d\}),E(\{a^{'},b^{'},c^{'},d^{'}\}))=2.
\end{equation}
By Claim 1, \eqref{equation Von(f) leq 7(2)} and \eqref{equation
cross(E0,E1)=2(2)},  there exists a region $h$ of $\ell\mbox{-}CQ_3$
with
\begin{equation}\label{equation CQ claim |Von(f)|=7}
|V_{on}(h)|=7.
\end{equation}
We may assume that $h$ is an unbounded region of $\ell\mbox{-}CQ_3$
(see Figure 3.4). By \eqref{equation CQ claim |Von(f)|=7}, at least
one edge of $E[\{a,b,c,d\},\{a^{'},b^{'},c^{'},d^{'}\}]$ would cross
some edge of $E(\{a,b,c,d\})\cup E(\{a^{'},b^{'},c^{'},d^{'}\})$.
This implies $\cross(E_{\ell})>2$, a contradiction with
\eqref{equation cross(E0,E1)=2(2)}. This proves Claim 2. \qed
\begin{figure}[h]
\centering
\includegraphics[scale=1.0]{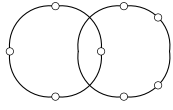}
\caption{\small{The graph for the proof of Claim 2}}
\end{figure}

\noindent\textbf{Claim 3.} {\sl There exists a region $f$ of
$\langle\{a,b,c,d\}\rangle$ with
$|V_{on}(f)|=|V_{in}(f;\langle\{a^{'},b^{'},c^{'},d^{'}\}\rangle)|=4$.}

{\sl Proof of Claim 3.} By Claim 2,
$\cross(E(\langle\{a,b,c,d\}\rangle),E(\langle\{a^{'},b^{'},c^{'},d^{'}\}\rangle))=0$.
This implies that all vertices of
$\langle\{a^{'},b^{'},c^{'},d^{'}\}\rangle$ lie in the same region
$f$ of $\langle\{a,b,c,d\}\rangle$, i.e.,
$|V_{in}(f;\langle\{a^{'},b^{'},c^{'},d^{'}\}\rangle)|=4$.

Now it remains to show $|V_{on}(f)|=4$. Suppose to the contrary that
$|V_{on}(f)|<4$, i.e., $\langle\{a,b,c,d\}\rangle$ crosses itself,
which implies
\begin{equation}\label{equation |Von(f)|=2}
|V_{on}(f)|=2.
\end{equation}
It follows that $\cross(E_{\ell})\geq
 \cross(E(\{a,b,c,d\}))+\cross(E(\{a,b,c,d\}),E[\{a,b,c,d\},\{a^{'},b^{'},c^{'},d^{'}\}])
\geq 1+2=3.$ By \eqref{equation 1 leq cr(El) leq 3(2)}, we have
$$\cross(E_{\ell})=3.$$ It follows from \eqref{equation |Von(f)|=2}
that $|V_{on}(h)|\leq 6$ for every region $h$ of $\ell\mbox{-}CQ_3$,
a contradiction with Claim 1. This proves Claim 3. \qed

By \eqref{equation CQ four 4-cycles}, \eqref{equation CQ four
4-cycles in particular 1} and \eqref{equation CQ four 4-cycles in
particular 2} and the structure of $CQ_3$, we conclude that
\begin{equation}\label{equation another two for four 4-cycles}
\mbox{ \ \ either \ \ }\{aa^{'},bb^{'},cd^{'},dc^{'}\}\subset
E_{\ell} \mbox{ \ \ or \ \ } \{ab^{'},ba^{'}, cc^{'},dd^{'}\}\subset
E_{\ell}.
\end{equation}

\noindent\textbf{Claim 4.} {\sl The edges
$ab,cd,a^{'}b^{'},c^{'}d^{'}$ are clean.}

{\sl Proof of Claim 4.} By \eqref{equation another two for four
4-cycles}, we may assume without loss of generality that
$$\{aa^{'},bb^{'},cd^{'},dc^{'}\}\subset E_{\ell}.$$ Since
$ab\in E(\{a,b,c,d\})\cap E( \{a,b,a^{'},b^{'}\})$, it follows from
Claim 2  that
\begin{equation}\label{equaiton CQ claim 4 ab no cross}
\cross(ab,E(\{a^{'},b^{'},c^{'},d^{'}\}))=\cross(ab,E(\{c,d,c^{'},d^{'}\}))=0.
\end{equation}
Because that no adjacent edges  cross each other, we have
\begin{equation}\label{equation CQ claim 4 ab no cross 2}
\cross(ab,\{ad,bc,aa^{'},bb^{'}\})=0.
\end{equation}
By \eqref{equaiton CQ claim 4 ab no cross} and \eqref{equation CQ
claim 4 ab no cross 2}, we conclude that $ab$ is clean.

A similar argument can be used to establish that the other three
edges $cd,a^{'}b^{'},c^{'}d^{'}$ are clean. This proves Claim 4.
\qed

By Claim 2, Claim 3 and Claim 4, there are three possible cases for
the drawing of $\langle\{a,b,c,d\}\rangle$ and
$\langle\{a^{'},b^{'},c^{'},d^{'}\}\rangle$, which are shown in
Figure 3.5.
\begin{figure}[h]
\centering
\includegraphics[scale=1.0]{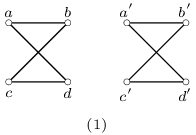} \hspace{50pt}
\includegraphics[scale=1.0]{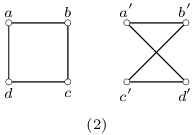} \hspace{50pt}
\includegraphics[scale=1.0]{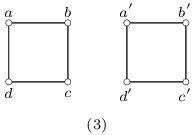}
\caption{\small{Possible cases for the drawing of
$\langle\{a,b,c,d\}\rangle$ and
$\langle\{a^{'},b^{'},c^{'},d^{'}\}\rangle$}}
\end{figure}

\noindent\textbf{Claim 5.} {\sl If $\{aa^{'},bb^{'}\}\subset
E_{\ell}$, then $\cross(\{aa^{'},bb^{'}\},E_{\ell}\setminus
\{aa^{'},bb^{'}\})+\cross(aa^{'},bb^{'})\geq 1$. If
$\{ab^{'},ba^{'}\}\subset E_{\ell}$, then
$\cross(\{ab^{'},ba^{'}\},E_{\ell}\setminus
\{ab^{'},ba^{'}\})+\cross(ab^{'},ba^{'})\neq 1$.}

{\sl Proof of Claim 5.} Suppose $\{aa^{'},bb^{'}\}\subset E_{\ell}$
and $\cross(\{aa^{'},bb^{'}\},E_{\ell}\setminus
\{aa^{'},bb^{'}\})+\cross(aa^{'},bb^{'})=0$. Then one edge of
$\{cd,c^{'}d^{'}\}$ must be drawn in the interior of $\langle
\{a,b,a^{'},b^{'}\}\rangle$ and the other in the exterior of
$\langle \{a,b,a^{'},b^{'}\}\rangle$, which implies that
$\cross(E(\{a,b,a^{'},b^{'}\}),E(\{c,d,c^{'},d^{'}\}))>0$, a
contradiction with Claim 2.

Suppose $\{ab^{'},ba^{'}\}\subset E_{\ell}$ and
$\cross(\{ab^{'},ba^{'}\},E_{\ell}\setminus
\{ab^{'},ba^{'}\})+\cross(ab^{'},ba^{'})=1$. We may assume $ab^{'}$
is not clean. By Claim 2, we have
$$\cross(ab^{'},\{cd,c^{'}d^{'}\})=0.$$ Because that no adjacent
edges cross each other, we have
$$\cross(ab^{'},\{ad,ab,a^{'}b^{'},b^{'}c^{'}\})=0.$$ It follows that
$$\cross(ab^{'},\{bc,a^{'}d^{'},ba^{'}\})=1,$$ and thus, one edge of
$\{cd,c^{'}d^{'}\}$ must be drawn in the interior of $\langle
\{a,b,a^{'},b^{'}\}\rangle$ and the other in the exterior of
$\langle \{a,b,a^{'},b^{'}\}\rangle$, which implies
$\cross(E(\{a,b,a^{'},b^{'}\}),E(\{c,d,c^{'},d^{'}\}))>0$, a
contradiction with Claim 2. This proves Claim 5. \qed

\noindent\textbf{Claim 6.} {\sl $\cross(E(\{a,b,c,d\}))=0$ or
$\cross(E(\{a^{'},b^{'},c^{'},d^{'}\}))=0$.}

{\sl Proof of Claim 6.} Assume to the contrary that
$$\cross(E(\{a,b,c,d\}))=\cross(E(\{a^{'},b^{'},c^{'},d^{'}\}))=1$$
 (see Figure
3.2(1)). By \eqref{equation another two for four 4-cycles}, we may
assume without loss of generality that
$$\{aa^{'},bb^{'},cd^{'},dc^{'}\}\subset E_{\ell}.$$
By  Claim 5, we have $\cross(\{aa^{'},bb^{'}\},E_{\ell}\setminus
\{aa^{'},bb^{'}\})+\cross(aa^{'},bb^{'})\geq 1$,  and so
$\cross(E_{\ell})\geq 3$. Combined with \eqref{equation 1 leq cr(El)
leq 3(2)}, we have
\begin{equation}\label{equaiton CQ claim 6 el=3}
\cross(E_{\ell})=3,
\end{equation}
which implies that both $cd^{'}$ and $dc^{'}$ are clean. Combined
with Claim 4, we have that all edges of $E(\{c,d,c^{'},d^{'}\})$ are
clean. By Lemma \ref{Lemma |Von(f)| leq 6}, we have that
$|V_{on}(f)|\leq 6$ for every region $f$ of $\ell\mbox{-}CQ_3$. By
Claim 1 and \eqref{equaiton CQ claim 6 el=3}, we derive a
contradiction.  This proves Claim 6. \qed

By Claim 6 and Claim 3, we may assume without loss of generality
that $$\cross(E(\{a,b,c,d\}))=0$$ and that
$\langle\{a^{'},b^{'},c^{'},d^{'}\}\rangle$ lies in the exterior of
$\langle\{a,b,c,d\}\rangle$.

Suppose $\cross(E(\{a^{'},b^{'},c^{'},d^{'}\}))=1$ (see Figure
3.5(2)).

If $\{aa^{'},bb^{'},cd^{'},dc^{'}\}\subset E_{\ell}$, by Claim 5, we
conclude that
$$\cross(\{aa^{'},bb^{'}\},E_{\ell}\setminus
\{aa^{'},bb^{'}\})+\cross(aa^{'},bb^{'})\geq 1,$$ and similarly,
$$\cross(\{cd^{'},dc^{'}\},E_{\ell}\setminus
\{cd^{'},dc^{'}\})+\cross(cd^{'},dc^{'})\geq 1.$$ It follows that
$\cross(E_{\ell})\geq 3$, and thus, by \eqref{equation 1 leq cr(El)
leq 3(2)}, $\cross(E_{\ell})=3.$ By Claim 1, we have that there
exists a region $f$ of $\ell\mbox{-}CQ_3$ with $|V_{on}(f)|=8$.
Without loss of generality, we may assume that $f$ is the unbounded
region of $\ell\mbox{-}CQ_3$. This implies $\cross(E_{\ell})>3$ (see
Figure 3.6(1)), a contradiction with \eqref{equation 1 leq cr(El)
leq 3(2)}.

If $\{ab^{'},ba^{'}, cc^{'},dd^{'}\}\subset E_{\ell}$, by Claim 5,
we conclude that
$$\cross(\{aa^{'},bb^{'}\},E_{\ell}\setminus
\{aa^{'},bb^{'}\})+\cross(aa^{'},bb^{'})\neq 1,$$ and similarly,
$$\cross(\{cd^{'},dc^{'}\},E_{\ell}\setminus
\{cd^{'},dc^{'}\})+\cross(cd^{'},dc^{'})\neq 1.$$ By \eqref{equation
1 leq cr(El) leq 3(2)}, we have that
$\cross(\{aa^{'},bb^{'}\},E_{\ell}\setminus
\{aa^{'},bb^{'}\})+\cross(aa^{'},bb^{'})=0$ or
$\cross(\{cd^{'},dc^{'}\}$, $E_{\ell}\setminus
\{cd^{'},dc^{'}\})+\cross(cd^{'},dc^{'})=0$, say
$$\cross(\{cd^{'},dc^{'}\},E_{\ell}\setminus
\{cd^{'},dc^{'}\})+\cross(cd^{'},dc^{'})=0.$$ That is, $cc^{'}$ and
$dd^{'}$ are clean. Combined with Claim 4, we have that all edges of
$E(\{c,d,c^{'},d^{'}\})$ are clean. By Lemma \ref{Lemma |Von(f)| leq
6}, we have that $|V_{on}(f)|\leq 6$ for every region $f$ of
$\ell\mbox{-}CQ_3$. By \eqref{equation 1 leq cr(El) leq 3(2)} and
Claim 1, we conclude that $\cross(E_{\ell})=1$. It follows that
$ab^{'}$ and $ba^{'}$ are also clean. By Claim 4, there is only one
possible drawing of $\ell\mbox{-}CQ_3$ as shown in Figure 3.6(2).
\begin{figure}[h]
\centering
\includegraphics[scale=1.0]{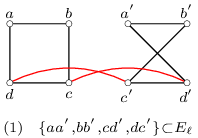} \hspace{50pt}
\includegraphics[scale=1.0]{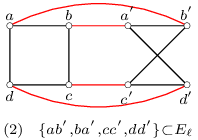}
\caption{\small{The cases for $\cross(ad,bc)=0$ and
$\cross(a^{'}d^{'},b^{'}c^{'})=1$}}
\end{figure}

Suppose $\cross(\langle\{a^{'},b^{'},c^{'},d^{'}\}\rangle)=0$  (see
Figure 3.5(3)). By \eqref{equation another two for four 4-cycles},
we may assume without loss of generality that
$$\{aa^{'},bb^{'},cd^{'},dc^{'}\}\subset E_{\ell}.$$ By
Claim 5, we have that $$\cross(\{aa^{'},bb^{'}\}, E_{\ell}\setminus
\{aa^{'},bb^{'}\})+\cross(aa^{'},bb^{'})\geq 1,$$ and similarly,
$$\cross(\{cd^{'},dc^{'}\},E_{\ell}\setminus
\{cd^{'},dc^{'}\})+\cross(cd^{'},dc^{'})\neq 1.$$

If $\cross(\{cd^{'},dc^{'}\},E_{\ell}\setminus
\{cd^{'},dc^{'}\})+\cross(cd^{'},dc^{'})\geq 2$, then
$\cross(E_{\ell})\geq 3$, and thus, by \eqref{equation 1 leq cr(El)
leq 3(2)}, $\cross(E_{\ell})=3.$ By Claim 1, we have that there
exists a region $f$ of $\ell\mbox{-}CQ_3$ with $|V_{on}(f)|=8$. We
may assume that $f$ is the unbounded region of $\ell\mbox{-}CQ_3$.
This implies $\cross(E_l)>3$ (see Figure 3.7(1)), a contradiction
with \eqref{equation 1 leq cr(El) leq 3(2)}.

If $\cross(\{cd^{'},dc^{'}\},E_{\ell}\setminus
\{cd^{'},dc^{'}\})+\cross(cd^{'},dc^{'})=0$, i.e., both $cd^{'}$ and
$dc^{'}$ are clean, by Claim 4, we have that all edges of
$E(\{c,d,c^{'},d^{'}\})$ are clean. By Lemma \ref{Lemma |Von(f)| leq
6}, we have that $|V_{on}(f)|\leq 6$ for every region $f$ of
$\ell\mbox{-}CQ_3$. By \eqref{equation 1 leq cr(El) leq 3(2)} and
Claim 1, we conclude that $\cross(E_{\ell})=1$. Combined with Claim
2 and Claim 4, by symmetry, there is only three possible drawings of
$\ell\mbox{-}CQ_3$ as shown in Figure 3.7(2)-(4).
\begin{figure}[h]
\centering
\includegraphics[scale=1.0]{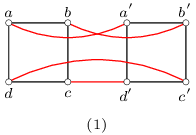} \hspace{5pt}
\includegraphics[scale=1.0]{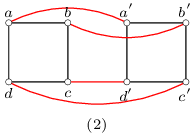} \hspace{5pt}
\includegraphics[scale=1.0]{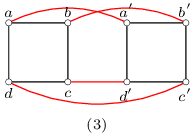} \hspace{5pt}
\includegraphics[scale=1.0]{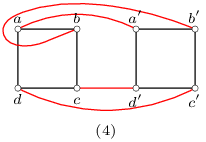}
\caption{\small{The cases for $\cross(ad,bc)=0$ and
$\cross(a^{'}d^{'},b^{'}c^{'})=0$}}
\end{figure}

Notice that the drawing shown in Figure 3.6(2) is isomorphic to the
drawing shown in Figure 3.7(3), and that the drawing shown in Figure
3.7(4) is isomorphic to the drawing shown in Figure 3.7(2). So we
need only to consider the drawings of $\ell\mbox{-}CQ_3$ shown in
Figure 3.7(2) and Figure 3.7(3).

From Figure 3.7(2) and 3.7(3), we see that
\begin{equation}\label{equation cr(El)=1(2)}
\cross(E_{\ell})=1
\end{equation}
and
\begin{equation}\label{equation max(Von(f))=5(2)}
|V_{on}(f)|\leq 5
\end{equation}
for every region $f$ of $\ell\mbox{-}CQ_3$. By \eqref{left leq
3.5(2)} and \eqref{equation cr(El)=1(2)}, we have that
\begin{equation}\label{equation cr(El,Elr) leq 2(2)}
\cross(E_{\ell},E_{\ell,r})\leq 2
\end{equation}
and
\begin{equation}\label{equation cr(El,Er)leq 5(2)}
\cross(E_{\ell},E_r)\leq 5.
\end{equation}
By Lemmas \ref{Lemma the number of emitting edges}, \ref{Lemma three
vertex set (6 1 1)}, \ref{Lemma four region geq 7} and \ref{Lemma
five regions geq 8}, we conclude that all vertices of
$r\mbox{-}CQ_3$ lie in at most three regions of $\ell\mbox{-}CQ_3$.

Suppose that all vertices of $r\mbox{-}CQ_3$ lie in the same region
of $\ell\mbox{-}CQ_3$. By \eqref{equation max(Von(f))=5(2)}, we have
$\cross(E_{\ell},E_{\ell,r})\geq 3$, a contradiction with
\eqref{equation cr(El,Elr) leq 2(2)}.

Suppose that all vertices of $r\mbox{-}CQ_3$ lie in exactly two
regions $f_1$ and $f_2$ of $\ell\mbox{-}CQ_3$. By Lemma \ref{Lemma
the number of emitting edges}, we have
$$\cross(E_{\ell},E_r)\geq 3.$$
It follows from \eqref{left leq 3.5(2)} and \eqref{equation
cr(El)=1(2)} that
\begin{equation}\label{equation CQ3 El,Elr leq 1}
\cross(E_{\ell},E_{\ell,r})\leq 1.
\end{equation} If
$\cross(E_{\ell},E_r)=3$, by Lemma \ref{Lemma the number of emitting
edges}, we have that
$\{|V_{in}(f_1,r\mbox{-}CQ_3)|,|V_{in}(f_2,r\mbox{-}CQ_3)|\}=\{1,7\}$.
It follows from \eqref{equation max(Von(f))=5(2)} that
$\cross(E_{\ell},E_{\ell,r})\geq 2$, a contradiction with
\eqref{equation CQ3 El,Elr leq 1}. Hence,
$$\cross(E_{\ell},E_r)\geq
4.$$ By \eqref{left leq 3.5(2)} and \eqref{equation cr(El)=1(2)}, we
have $\cross(E_{\ell},E_{\ell,r})=0$, which implies $V_{on}(f_1)\cup
V_{on}(f_2)=V(\ell\mbox{-}CQ_3)$. From the two drawings of Figure
3.7(2)-(3), we observe that the two regions containing all the
vertices of $\ell\mbox{-}CQ_3$ do not have common boundary. By Lemma
\ref{Lemma the number of emitting edges}, we have
$\cross(E_{\ell},E_r)\geq 8$, a contradiction with \eqref{equation
cr(El,Er)leq 5(2)}.

Suppose that all vertices of $r\mbox{-}CQ_3$ lie in exactly three
regions $f_1$, $f_2$ and $f_3$ of $\ell\mbox{-}CQ_3$. By \eqref{left
leq 3.5(2)}, \eqref{equation cr(El)=1(2)}, \eqref{equation
cr(El,Er)leq 5(2)} and Lemma \ref{Lemma three vertex set (6 1 1)},
we conclude that $\cross(E_{\ell},E_r)=5$,
$\cross(E_{\ell},E_{\ell,r})=0$ and
$\{|V_{in}(f_1,r\mbox{-}CQ_3)|,|V_{in}(f_2,r\mbox{-}CQ_3)|,|V_{in}(f_3,r\mbox{-}CQ_3)|\}=\{1,1,6\}$.
It follows from \eqref{equation max(Von(f))=5(2)} that
$\cross(E_{\ell},E_{\ell,r})\geq 1$, a contradiction with
$\cross(E_{\ell},E_{\ell,r})=0$.

This completes the proof of Theorem \ref{theorem cr(CQ4)=8}.
\end{proof}

\section{Crossing number of $LTQ_4$}

\indent \indent In Figure 4.1, we give a drawing of $LTQ_4$ with 10
crossings, Hence, we have the following

\begin{lemma} \label{lemma upper bound of LTQ4}
$cr(LTQ_4)\leq 10$.
\end{lemma}
\begin{figure}[ht]
\centering
\includegraphics[scale=0.9]{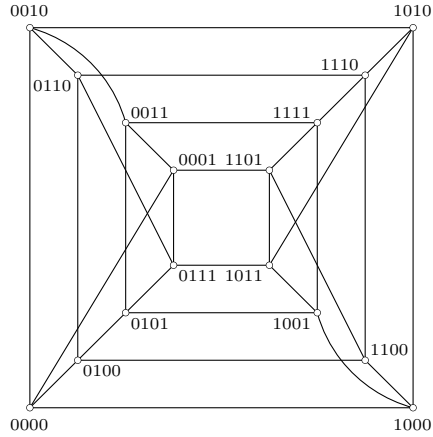}
\caption{\small{A good drawing of $LTQ_4$ with 10 crossings}}
\end{figure}

In the rest of this section, we shall prove that the  value of
$cr(LTQ_4)$ is exactly equal to 10. We rename the vertices of
$LTQ_4$ as shown in Figure 4.2. Let $\ell\mbox{-}LTQ_3$ and
$r\mbox{-}LTQ_3$ be the subgraphs induced by vertex subset
$\{v_i:0\leq i\leq 7\}$ and by vertex subset $\{v_i:8\leq i\leq
15\}$, respectively. Note that both $\ell\mbox{-}LTQ_3$ and
$r\mbox{-}LTQ_3$ are isomorphic to $LTQ_3$.
For convenience, we abbreviate
$$\begin {array}{llll}
E_{\ell}&=&E(\ell\mbox{-}LTQ_3),\\
E_r&=&E(r\mbox{-}LTQ_3),\\
E_{\ell,r}&=&E[V(\ell\mbox{-}LTQ_3),V(r\mbox{-}LTQ_3)],\\
\end {array}$$
\begin{figure}[ht]
\centering
\includegraphics[scale=0.9]{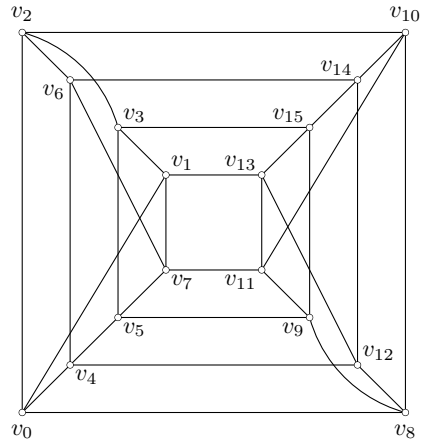}
\caption{\small{ $LTQ_4$ with the renamed vertices}}
\end{figure}

By Lemma \ref{lemma counting the number of crossings}, we have the
following

\begin{lemma}\label{Lemma LTQ4 total crossing and left  }
Let $D$ be a good drawing
of $LTQ_4$. Then
$$\nu(D)=\cross(E_{\ell})+\cross(E_r)+\cross(E_{\ell,r})+\cross(E_{\ell},E_r)+\cross(E_{\ell},E_{\ell,r})+\cross(E_r,E_{\ell,r}).$$
\end{lemma}

Before proceeding, we need some preliminaries.

\begin{definition}
Let $\pi:V(\ell\mbox{-}LTQ_3)\rightarrow V(r\mbox{-}LTQ_3)$ be a
bijection defined as follows: $u\rightarrow \pi(u)$ where $\pi(u)\in
V(r\mbox{-}LTQ_3)$ is adjacent to $u$.
\end{definition}

\begin{observation} \label{observation 4-cycle match}
Let $X$ be a subset of $V(LTQ_3)$. The following
conditions hold. \\
\noindent(i) If $\langle X\rangle \cong P_2$, then there exists a
subset $Y\subseteq V(LTQ_3)\setminus X$ such that $\langle X\cup
Y\rangle
\cong C_4$. \\
\noindent(ii) If $\langle X\rangle \cong C_4$, then $\langle
V(LTQ_3)\setminus X\rangle \cong C_4$ and there exists a matching of
cardinality four between $X$ and $V(LTQ_3)\setminus X$.
\end{observation}

\begin{observation} \label{observation vertex partition of LTQ4} Let
$u_1,u_2,u_3,u_4,u_5\in V(\ell\mbox{-}LTQ_3)$. The following
conditions
hold. \\
(i) If $\langle\{u_1,u_2,u_3,u_4,u_5\}\rangle\cong C_5$, then
$\langle \{\pi(u_1),\pi(u_2),\pi(u_3),\pi(u_4),\pi(u_5)\}\setminus
\{\pi(u_i)\}\rangle\not\cong C_4$ for any $1\leq i\leq 5$. \\
(ii) If $\langle \{u_1,u_2,u_3,u_4\}\rangle\cong C_4$, then
$\langle\{\pi(u_1),\pi(u_2),\pi(u_3),\pi(u_4),\pi(u_5)\}\setminus
\{\pi(u_i)\}\rangle \not \cong C_4$ for any $1\leq i\leq 4$.
\end{observation}

\begin{observation} \label{observation the structure of P3 in l-LTQ3}
Let $u_1,u_2,u_3\in V(\ell\mbox{-}LTQ_3)$ with $u_1u_2,u_2u_3\in
E(LTQ_4)$. If $\pi(u_1)\pi(u_2)$, $\pi(u_2)\pi(u_3)\in E(LTQ_4)$,
then there exists $u\in V(\ell\mbox{-}LTQ_3)$ such that $\langle
\{u_1,u_2,u_3,u\}\rangle\cong C_4$.
\end{observation}

\begin{observation} \label{observation the paths P4 of LTQ_4}
Let $u_1,u_2,u_3,u_4\in V(\ell\mbox{-}LTQ_3)$ with $u_1u_2,u_3u_4\in
E(\ell\mbox{-}LTQ_3)$. Assume that, for any $V\subseteq
V(\ell\mbox{-}LTQ_3)$ such that $\langle V\rangle\cong C_4$, either
$V\cap \{u_1,u_2,u_3,u_4\}=\{u_1,u_2\}$ or $V\cap
\{u_1,u_2,u_3,u_4\}=\{u_3,u_4\}$. Then there exist two
vertex-disjoint 4-paths $P(u_{i_1},u_{j_1}),P(u_{i_2},u_{j_2})$ of
$LTQ_4$ such that $P(u_{i_1},u_{j_1})\cap
P(u_{i_2},u_{j_2})=\emptyset$ where $\{i_1,i_2\}=\{1,2\}$ and
$\{j_1,j_2\}=\{3,4\}$ (see Figure 4.3).
\end{observation}
\begin{figure}[ht]
\centering
\includegraphics[scale=1.0]{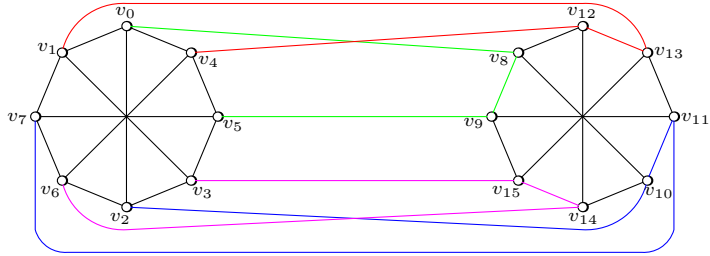}
\caption{\small{Four vertex-disjoint 4-paths}}
\end{figure}

\begin{observation} \label{observation the relation of paths}
Let $C$ be a closed curve, and let $u_1,u_2,u_3,u_4$ be points drawn
on $C$ along clockwise direction. Denote $L_{i,j}$ to be a line
drawn between $u_i$ and $u_j$ in the interior of $C$, where $1\leq
i<j\leq 4$. Then the number of crossings produced by $L_{1,2}$ and
$L_{3,4}$ is even, and the number of crossings produced by $L_{1,3}$
and $L_{2,4}$ is odd.
\end{observation}

\begin{lemma} \label{Lemma crossings above 4-cycle and 5-cycle}
Let $C$ be a closed curve, and let $u_1,u_2,u_3,u_4,u_5$ be vertices
in the interior of $C$ with
$\langle\{u_1,u_2,u_3,u_4,u_5\}\rangle\cong C_5$. For every vertex
$u_i$ there exists two edges $e_i^1,e_i^2$ incident with $u_i$ and
vertices $v_i^1, v_i^2$ which are lying on the curve $C$ or in the
exterior of $C$. Let $u_i^{'}$ be the vertex $v_i^1$ or the crossing
produced by $e_i^1$ and $C$, according to $v_i^1$ lying on $C$ or
not. Let $\mathscr{C}_1,\mathscr{C}_2,\ldots, \mathscr{C}_5$ be the
line segments of $C$ divided by $u_i^{'}$. Assume that there exists
$i\in \{1,2,\ldots,5\}$ such that every edge of
$\{e_1^2,e_2^2,\ldots,e_5^2\}$ crosses $\mathscr{C}_i$, and that
$\{e_1^1,e_2^1,\ldots,e_5^1\}$ do not cross each other. Then the
number of crossings produced among the edges of
$E(\langle\{u_1,u_2,u_3,u_4,u_5\}\rangle)\cup
\{e_1^1,e_2^1,e_3^1,e_4^1,e_5^1\}\cup
\{e_1^2,e_2^2,e_3^2,e_4^2,e_5^2\}$ is at least 3.
\end{lemma}

\begin{proof} For convenience, let $E_c=E(\langle\{u_1,u_2,u_3,u_4,u_5\}\rangle)$,
$E_1=\{e_1^1,e_2^1,e_3^1,e_4^1,e_5^1\}$ and
$E_2=\{e_1^2,e_2^2,e_3^2,e_4^2,e_5^2\}$. Let $\beta$ be the number
of crossings produced among the edges of $E_c\cup E_1\cup E_2$.
Assume without loss of generality that $u_1,u_2,u_3,u_4,u_5$ are
arranged along count-clockwise direction and that every edge of
$E_2$ crosses the line segment between $u_1^{'}$ and $u_2^{'}$.

Suppose to the contrary that $\beta\leq 2$. If $u_1u_2\in E_c$, then
the path consist of $e_1^1,e_2^1,u_1u_2$ has to cross one of
$\{e_i^1,e_i^2\}$, where $i=3,4,5$. This implies that $\beta\geq 3$
(see Figure 4.4(1)), a contradiction. Hence, we have
$$u_1u_2\not\in E_c.$$ If $u_1u_3\in E_c$ or $u_2u_5\in E_c$, say
$u_1u_3\in E_c$, then $u_2u_4\in E_c$ or $u_2u_5\in E_c$. It follows
that the path consist of $e_1^1,e_3^1,u_1u_3$ has to cross one of
$\{e_i^1,e_i^2\}$ where $i=4,5$, and the path consist of
$e_2^1,e_j^1,u_2u_j$ has to cross one of $\{e_1^1,e_3^1,u_1u_3\}$
where $j=4$ or $j=5$. This implies that $\beta\geq 3$ (see Figure
4.4(2)), a contradiction. Hence, we have
$$u_1u_3,u_2u_5\not\in E_c.$$ It follows that $u_1u_4,u_1u_5,u_2u_3,u_2u_4,u_3u_5\in
E_c$. Thus, we have that the path consist of $e_1^1,e_4^1,u_1u_4$
has to cross one of $\{e_5^1,e_5^2\}$, the path consist of
$e_2^1,e_4^1,u_2u_4$ has to cross one of $\{e_3^1,e_3^2\}$, and the
path consist of $e_3^1,e_5^1,u_3u_5$ has to cross one of
$\{e_4^1,e_4^2\}$. This implies that $\beta\geq 3$ (see Figure
4.4(3)), a contradiction.
\begin{figure}[h]
\centering
\includegraphics[scale=1.0]{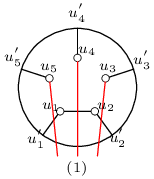} \hspace{20pt}
\includegraphics[scale=1.0]{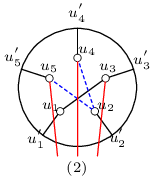} \hspace{20pt}
\includegraphics[scale=1.0]{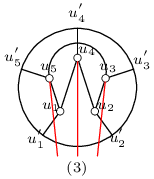}
\caption{\small{The graphs for the proof of Lemma \ref{Lemma
crossings above 4-cycle and 5-cycle}}}
\end{figure}
\end{proof}

\begin{theorem} \label{theorem cr(LTQ4)=10}
$cr(LTQ_4)=10$.
\end{theorem}

\begin{proof} By Lemma \ref{lemma upper bound of LTQ4}, it suffices to prove that
$cr(LTQ_4)\geq 10$. Suppose to the contrary that
\begin{equation}\label{equation cr(LTQ4) leq 9}
cr(LTQ_4)\leq 9,
\end{equation}
i.e., there exists a good drawing $D$ of $LTQ_4$ such that
$\nu(D)\leq 9$. Without loss of generality, we may assume
$\cross(E_{\ell})+\cross(E_{\ell},E_{\ell,r})\leq
\cross(E_r)+\cross(E_r,E_{\ell,r})$.
Since $\nu(D)\leq 9$, by Lemma
\ref{Lemma LTQ4 total crossing and left }, we have
\begin{equation}\label{left leq 3.5}
\cross(E_{\ell})+\cross(E_{\ell},E_{\ell,r})+0.5\cross(E_{\ell,r})+0.5\cross(E_{\ell},E_r)\leq
4.5
\end{equation}
By Observation \ref{Observation G3 contains K3,3}, we have
$\cross(E_{\ell})\geq 1$ and $\cross(E_r)\geq 1$. It follows from
\eqref{left leq 3.5} that
\begin{equation}\label{equation 1 leq cr(El) leq 4}
1\leq \cross(E_{\ell})\leq 4.
\end{equation}
By \eqref{left leq 3.5} and \eqref{equation 1 leq cr(El) leq 4}, we
have
\begin{equation}\label{equation cr(El,Elr) leq 3}
\cross(E_{\ell},E_{\ell,r})\leq 3.
\end{equation}

By Observation \ref{observation vertex partition of LTQ3}, we may
assume without loss of generality that there exists a permutation
$a,b,c,d,a^{'},b^{'},c^{'},d^{'}$ of $V(\ell\mbox{-}LTQ_3)$
satisfying that
\begin{equation}\label{equation LTQ four 4-cycles}
\langle\{a,b,c,d\}\rangle\cong
\langle\{a^{'},b^{'},c^{'},d^{'}\}\rangle\cong \langle
\{a,b,a^{'},b^{'}\}\rangle\cong \langle
\{c,d,c^{'},d^{'}\}\rangle\cong C_4,
\end{equation}
 and moreover,
\begin{equation}\label{equation LTQ four 4-cycles in particular 1}
\langle\{a,b,c,d\}\rangle=abcda,
\end{equation}
\begin{equation}\label{equation LTQ four 4-cycles in particular 2}
\langle\{a^{'},b^{'},c^{'},d^{'}\}\rangle=a^{'}b^{'}c^{'}d^{'}a^{'}.
\end{equation}

By \eqref{equation LTQ four 4-cycles}, \eqref{equation LTQ four
4-cycles in particular 1} and \eqref{equation LTQ four 4-cycles in
particular 2} and the structure of $LTQ_3$, we conclude that
\begin{equation}\label{equation LTQ another two for four 4-cycles}
\mbox{ \ \ either \ \ }\{aa^{'},bb^{'},cd^{'},dc^{'}\}\subset
E_{\ell} \mbox{ \ \ or \ \ } \{ab^{'},ba^{'}, cc^{'},dd^{'}\}\subset
E_{\ell}.
\end{equation}

\begin{claim} \label{claim cr(El) geq 3} {\sl If
$\cross(E_{\ell})=4$ then there exists a region $f$ of
$\ell\mbox{-}LTQ_3$ with $|V_{on}(f)|=8$, and if
$\cross(E_{\ell})=3$ then there exists a region $f$ of
$\ell\mbox{-}LTQ_3$ with $|V_{on}(f)|\geq 7$.}
\end{claim}

{\sl Proof of Claim \ref{claim cr(El) geq 3}.} Suppose
$\cross(E_{\ell})=4$. By \eqref{left leq 3.5} and Lemma \ref{Lemma
the number of emitting edges}, we conclude that
$\cross(E_{\ell},E_{\ell,r})=0$ and $\cross(E_{\ell},E_r)=0$. It
follows that there exists a region $f$ of $\ell\mbox{-}LTQ_3$ with
$|V_{on}(f)|=8$.

Suppose that $\cross(E_{\ell})=3$ and
\begin{equation}\label{equation V_on(f) leq 6}
|V_{on}(f)|\leq 6
\end{equation}
for every region $f$ of $\ell\mbox{-}LTQ_3$. By \eqref{left leq
3.5}, we have
$$\cross(E_{\ell},E_{\ell,r})+0.5\cross(E_{\ell},E_r)\leq 1.5.$$ By
\eqref{equation V_on(f) leq 6} and Lemma \ref{Lemma the number of
emitting edges}, we conclude that $\cross(E_{\ell},E_r)=3$ and there
exists a region $h$ of $\ell\mbox{-}LTQ_3$ such that
$|V_{in}(h,\ell\mbox{-}LTQ_3)|=7$, which is a contradiction with
\eqref{equation V_on(f) leq 6}. This proves Claim \ref{claim cr(El)
geq 3}. \qed

\begin{claim} \label{claim two cycles can not cross}
{\sl
$\cross(E(\{a,b,c,d\}),E(\{a^{'},b^{'},c^{'},d^{'}\}))=\cross(E(\{a,b,a^{'},b^{'}\}),E(\{c,d,
c^{'},d^{'}\}))=0$.}
\end{claim}

{\sl Proof of Claim \ref{claim two cycles can not cross}.} Without
loss of generality, we may assume to the contrary that
\begin{equation}\label{equation LTQ Claim 2 assume contrary}
\cross(E(\{a,b,c,d\}),E(\{a^{'},b^{'},c^{'},d^{'}\}))\neq 0.
\end{equation}
By \eqref{equation 1 leq cr(El) leq 4}, \eqref{equation LTQ Claim 2
assume contrary} and Observation \ref{Observation two
vertex-disjoint cycles}, we have
$$\cross(E(\{a,b,c,d\}),E(\{a^{'},b^{'},c^{'},d^{'}\}))\in\{2,4\}.$$
\indent Suppose
$\cross(E(\{a,b,c,d\}),E(\{a^{'},b^{'},c^{'},d^{'}\}))=4$. Then
$\cross(E_{\ell})=4$. By Claim \ref{claim cr(El) geq 3}, we have
that all the vertices of $\ell\mbox{-}LTQ_3$ are on the boundary of
the same region. Therefore, there exist only three possible drawings
as shown in Figure 4.5. It follows that either $cr(E_{\ell})\geq 5$
or the graph satisfying $\cross(E_{\ell})=4$ is isomorphic to $Q_3$,
a contradiction. Hence,
\begin{equation}\label{equation cross(E0,E1)=2}
\cross(E_{\ell})\geq
\cross(E(\{a,b,c,d\}),E(\{a^{'},b^{'},c^{'},d^{'}\}))=2.
\end{equation}
By \eqref{left leq 3.5} and \eqref{equation cross(E0,E1)=2}, we have
\begin{equation}\label{equation cr(El,Elr)+0.5cr(El,Er)leq 2.5}
\cross(E_{\ell},E_{\ell,r})+0.5\cross(E_{\ell},E_r)\leq 2.5.
\end{equation}

\begin{figure}[ht]
\centering
\includegraphics[scale=1.0]{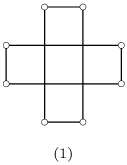} \hspace{20pt}
\includegraphics[scale=1.0]{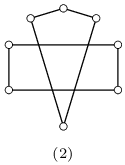} \hspace{20pt}
\includegraphics[scale=1.0]{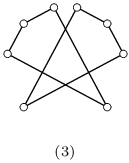}
\caption{\small{Three possible drawings of $\ell\mbox{-}LTQ_3$ for
$\cross(E_{\ell})=4$}}
\end{figure}

Let $f_{max}$ be the region of $\ell\mbox{-}LTQ_3$ with
$|V_{on}(f_{max})|=\max\limits_f\{|V_{on}(f)|\}$, where $f$ runs
over all regions of $\ell\mbox{-}LTQ_3$. Without loss of generality,
we may assume that $f_{max}$ is the the unbounded region of
$\ell\mbox{-}LTQ_3$. By \eqref{equation LTQ Claim 2 assume
contrary}, we have
$$1\leq |V_{on}(f_{max})|\leq 7.$$
\indent Suppose $|V_{on}(f_{max})|=7$. Then at least one edge of
$E[\{a,b,c,d\},\{a^{'},b^{'},c^{'},d^{'}\}]$ would cross some edge
of $E(\{a,b,c,d\})\cup E(\{a^{'},b^{'},c^{'},d^{'}\})$. This implies
$\cross(E_{\ell})\geq 4$. By \eqref{equation 1 leq cr(El) leq 4}, we
have $\cross(E_{\ell})=4$. By Claim \ref{claim cr(El) geq 3}, we
have that there exists a region $f$ of $\ell\mbox{-}LTQ_3$ with
$|V_{on}(f)|=8$, a contradiction with $|V_{on}(f_{max})|=7$. Hence,
$$|V_{on}(f_{max})|\leq 6.$$ By Claim \ref{claim cr(El) geq 3}, we
conclude that $\cross(E_{\ell})\leq 2$. It follows from
\eqref{equation cross(E0,E1)=2} that
\begin{equation}\label{equation cr(El)=2}
\cross(E_{\ell})=2
\end{equation}
and all the edges of
$E^{'}=E[\{a,b,c,d\},\{a^{'},b^{'},c^{'},d^{'}\}]$ are clean, that
is,
\begin{equation}\label{equation cr(E0)= cdots =cr(E1,E01)=0}
\cross(E^{'})+\cross(E^{'},E_{\ell}\setminus E^{'})=0.
\end{equation}
\indent Let $n_1$ be the number of vertices of
$\langle\{a^{'},b^{'},c^{'},d^{'}\rangle$ lying in the inner space
of $\langle\{a,b,c,d\}\rangle$, and let $n_2$ be the number of
vertices of $\langle\{a,b,c,d\}\rangle$ lying in the inner space of
$\langle\{a^{'},b^{'},c^{'},d^{'}\rangle$.

Suppose $|V_{on}(f_{max})|=6$. Then $(n_1,n_2)\in
\{(0,1),(1,0),(1,1),(0,2),(2,0)\}$.
\begin{figure}[h]
\centering
\includegraphics[scale=1.0]{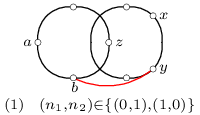} \hspace{40pt}
\includegraphics[scale=1.0]{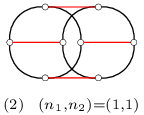} \hspace{40pt}
\includegraphics[scale=1.0]{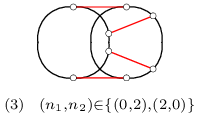}
\caption{\small{Three cases for $|V_{on}(f_{max})|=6$}}
\end{figure}

For $(n_1,n_2)\in \{(0,1),(1,0)\}$ (see Figure 4.6(1)), we have that
at least one of the vertices $x$ and $y$, say vertex $y$, is not
adjacent to vertex $z$. Then vertex $y$ must be adjacent to vertex
$b$. It follows that $\cross(E^{'})+\cross(E^{'},E_{\ell}\setminus
E^{'})\geq 1$, a contradiction with \eqref{equation cr(E0)= cdots
=cr(E1,E01)=0}.

For $(n_1,n_2)\in\{(1,1),(0,2),(2,0)\}$ (see Figure 4.6(2)-(3)), by
\eqref{equation cr(El)=2} and \eqref{equation cr(E0)= cdots
=cr(E1,E01)=0}, the edges of
$E[\{a,b,c,d\},\{a^{'},b^{'},c^{'},d^{'}\}]$ must be drawn as shown
in Figure 4.4(2)-(3). It is easy to see that the graphs shown in
Figure 4.4(2)-(3) are isomorphic to $Q_3$, a contradiction. Hence,
\begin{equation}\label{equation V_on(fmax) leq 5}
1\leq |V_{on}(f_{max})|\leq 5.
\end{equation}
\indent By \eqref{equation cr(El,Elr)+0.5cr(El,Er)leq 2.5}, we have
$\cross(E_{\ell},E_{\ell,r})\leq 2$. If $\cross(E_{\ell},E_r)=0$, by
\eqref{equation V_on(fmax) leq 5}, we have
$\cross(E_{\ell},E_{\ell,r})\geq 3$, a contradiction with
$\cross(E_{\ell},E_{\ell,r})\leq 2$. Hence,
$\cross(E_{\ell},E_r)\neq 0$. By \eqref{equation
cr(El,Elr)+0.5cr(El,Er)leq 2.5}, \eqref{equation V_on(fmax) leq 5},
Lemma \ref{Lemma the number of emitting edges}, Lemma \ref{Lemma
three vertex set (6 1 1)}, Lemma \ref{Lemma four region geq 7} and
Lemma \ref{Lemma five regions geq 8}, we conclude that
\begin{equation}\label{equation 4 leq cr(El,Er)
leq 5} 4\leq \cross(E_{\ell},E_r)\leq 5,
\end{equation}
and there exist exactly two regions $f_1,f_2$ of $\ell\mbox{-}LTQ_3$
with $\mathcal {B}(f_1,f_2)=1$ such that
$V(r\mbox{-}LTQ_3)=V_{in}(f_1;r\mbox{-}LTQ_3)\cup
V_{in}(f_2;r\mbox{-}LTQ_3)$ and
\begin{equation}\label{equation Von(f1)cupVon(f2)=V}
V(\ell\mbox{-}LTQ_3)=V_{on}(f_1)\cup V_{on}(f_2).
\end{equation}
We may assume without loss of generality that $$|V_{on}(f_1)|\geq
|V_{on}(f_2)|.$$ Combined with \eqref{equation V_on(fmax) leq 5}, we
have
$$|V_{on}(f_1)|=|V_{on}(f_{max})|\in\{4,5\}.$$
\indent Suppose $|V_{on}(f_1)|=4$. By \eqref{equation
Von(f1)cupVon(f2)=V}, we have that
\begin{equation}\label{equation Von(f1)cap Von(f2)=empty}
V_{on}(f_1)\cap V_{on}(f_2)=\emptyset.
\end{equation}
Let $x$ and $y$ be the two points at which
$\langle\{a,b,c,d\}\rangle$ and $E(\{a^{'},b^{'},c^{'},d^{'}\})$
cross each other. Let $\mathcal{A}_1,\mathcal{A}_3$ be the two parts
of 4-cycle $\langle\{a,b,c,d\}\rangle$ divided by $x$ and $y$, and
let $\mathcal{A}_2,\mathcal{A}_4$ be the two parts of 4-cycle
$\langle\{a^{'},b^{'},c^{'},d^{'}\}\rangle$ divided by $x$ and $y$
(see Figure 4.7).
\begin{figure}[h]
\centering
\includegraphics[scale=1.0]{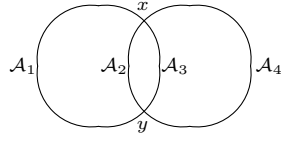}
\caption{\small{The graph for the proof of $|V_{on}(f_1)|=4$}}
\end{figure}
Let $\mathcal{M}_i$ be the vertex set that consists of all vertices
lying on $\mathcal{A}_i$ for $i=1,2,3,4$. Let $C$ be the common
boundary of the regions $f_1$ and $f_2$.  Note that
$$\mathcal{M}_1\cup \mathcal{M}_3=\{a,b,c,d\}$$ and $$\mathcal{M}_2\cup
\mathcal{M}_4=\{a^{'},b^{'},c^{'},d^{'}\}.$$ By \eqref{equation
Von(f1)cap Von(f2)=empty},  $C$ contains no vertex of
$\ell\mbox{-}LTQ_3$. Combined with \eqref{equation cross(E0,E1)=2}
and \eqref{equation cr(El)=2}, we conclude that
 $C$ is eaxctly one of  $\mathcal{A}_1,\mathcal{A}_2,\mathcal{A}_3$ and
 $\mathcal{A}_4$.  If $C$ is  $\mathcal{A}_1$ or
$\mathcal{A}_3$, we have $V_{on}(f_1)\cup V_{on}(f_2)\subseteq
  \mathcal{M}_2\cup\mathcal{M}_4$, a contradiction with \eqref{equation Von(f1)cupVon(f2)=V}.
If $C$ is  $\mathcal{A}_2$ or $\mathcal{A}_4$, we have
$V_{on}(f_1)\cup V_{on}(f_2)\subseteq
  \mathcal{M}_1\cup\mathcal{M}_3$, a contradiction with \eqref{equation Von(f1)cupVon(f2)=V}.
Therefore,
\begin{equation}\label{equation V_on(f2) leq V_on(f1)}
|V_{on}(f_1)|=5.
\end{equation}
Then $(n_1,n_2)\in
\{(0,1),(1,0),(1,1),(0,2),(2,0),(1,2),(2,1),(0,3),(3,0)\}$. Applying
\eqref{equation cr(El)=2} and \eqref{equation cr(E0)= cdots
=cr(E1,E01)=0}, by symmetry, the edges of
$E[\{a,b,c,d\},\{a^{'},b^{'},c^{'},d^{'}\}]$ must be drawn as Figure
4.8. Since $\mathcal {B}(f_1,f_2)=1$, then the only possible drawing
of $\ell\mbox{-}LTQ_3$ satisfying \eqref{equation
Von(f1)cupVon(f2)=V} and \eqref{equation V_on(f2) leq V_on(f1)} is
shown in Figure 4.8(3).
\begin{figure}[h]
\centering
\includegraphics[scale=1.0]{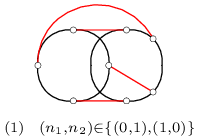} \hspace{30pt}
\includegraphics[scale=1.0]{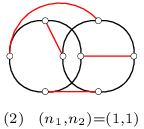} \hspace{30pt}
\includegraphics[scale=1.0]{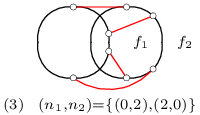}

\vspace{10pt}

\includegraphics[scale=1.0]{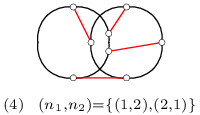} \hspace{30pt}
\includegraphics[scale=1.0]{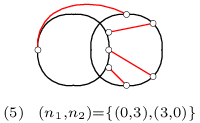}
\caption{\small{Five cases for $|V_{on}(f_{max})|=5$}}
\end{figure}

From Figure 4.8(3), we see that
$(|V_{in}(f_1;r\mbox{-}LTQ_3)|,|V_{in}(f_2;r\mbox{-}LTQ_3)|)\in\{(4,4),(3,5),(5,3)\}$,
$\langle V_{on}(f_1)\rangle \cong \langle V_{on}(f_2)\rangle \cong
C_5$ and the common boundary of $f_1$ and $f_2$ is an edge.

For
$(|V_{in}(f_1;r\mbox{-}LTQ_3)|,|V_{in}(f_2;r\mbox{-}LTQ_3)|)=(4,4)$,
by Observation \ref{observation vertex partition of LTQ4}(i), we
have that $\langle V_{in}(f_1;r\mbox{-}LTQ_3)\rangle\not\cong C_4$.
It follows from Lemma \ref{Lemma the number of emitting edges} that
$\cross(E_{\ell},E_r)\geq 6$, a contradiction with \eqref{equation 4
leq cr(El,Er) leq 5}.

For
$(|V_{in}(f_1;r\mbox{-}LTQ_3)|,|V_{in}(f_2;r\mbox{-}LTQ_3)|)\in\{(3,5),(5,3)\}$,
say $(|V_{in}(f_1;r\mbox{-}LTQ_3)|,|V_{in}$
$(f_2;r\mbox{-}LTQ_3)|)=(3,5)$, by \eqref{equation 4 leq cr(El,Er)
leq 5} and Lemma \ref{Lemma the number of emitting edges}, we have
$$\cross(E_{\ell},E_r)=5.$$
Since $V_{on}(f_2)\cong C_5$, it follows from Lemma \ref{Lemma the
number of emitting edges} that $\langle
V_{in}(f_2;r\mbox{-}LTQ_3)\rangle\cong C_5$ or $\langle
V_{in}(f_2;r\mbox{-}LTQ_3)\setminus \{v\}\rangle\cong C_4$ for some
$v\in V_{in}(f_2;r\mbox{-}LTQ_3)$. By Observation \ref{observation
vertex partition of LTQ4}, we conclude that $\langle
V_{in}(f_2;r\mbox{-}LTQ_3)\rangle\cong C_5$. By \eqref{equation
cr(LTQ4) leq 9}, \eqref{left leq 3.5}, \eqref{equation cr(El)=2} and
Lemma \ref{Lemma LTQ4 total crossing and left }(i), we have that
$\cross(E_{\ell,r})=0$,
$$\cross(E_r)+\cross(E_r,E_{\ell,r})\leq 2$$ and the common boundary is crossed
exactly five times by the edges of $E_r$. It follows from Lemma
\ref{Lemma crossings above 4-cycle and 5-cycle} that
$\cross(E_r)+\cross(E_r,E_{\ell,r})\geq 3$, a contradiction. This
proves Claim \ref{claim two cycles can not cross}. \qed

Similar arguments can be used to establish the following Claims 3-6.

\begin{claim} \label{claim |Von(f)|=|Vin(f;C)|=4}
{\sl There exists a region $f$ of $\langle\{a,b,c,d\}\rangle$ with
$|V_{on}(f)|=|V_{in}(f;\langle\{a^{'},b^{'},c^{'},d^{'}\}\rangle)|=4$.}
\end{claim}

\begin{claim} \label{claim four edges are clean}
{\sl The edges $ab,cd,a^{'}b^{'},c^{'}d^{'}$ are clean.}
\end{claim}

\begin{claim} \label{claim two edges cross geq 1 or neq 1}
{\sl If $\{aa^{'},bb^{'}\}\subset E_{\ell}$, then
$\cross(\{aa^{'},bb^{'}\},E_{\ell}\setminus
\{aa^{'},bb^{'}\})+\cross(aa^{'},bb^{'})\geq 1$. If
$\{ab^{'},ba^{'}\}\subset E_{\ell}$, then
$\cross(\{ab^{'},ba^{'}\},E_{\ell}\setminus
\{ab^{'},ba^{'}\})+\cross(ab^{'},ba^{'})\neq 1$.}
\end{claim}

\begin{claim} \label{claim two cycles do not cross itself}
{\sl $\cross(E(\{a,b,c,d\}))=0$ or
$\cross(E(\{a^{'},b^{'},c^{'},d^{'}\}))=0$.}
\end{claim}

By Claim \ref{claim two cycles do not cross itself} and Claim
\ref{claim |Von(f)|=|Vin(f;C)|=4}, we may assume without loss of
generality that $$\cross(E(\{a,b,c,d\}))=0$$ and that
$\langle\{a^{'},b^{'},c^{'},d^{'}\}\rangle$ lies in the exterior of
$\langle\{a,b,c,d\}\rangle$.

Suppose $\cross(\langle\{a^{'},b^{'},c^{'},d^{'}\}\rangle)=1$. By an
argument similar as Theorem \ref{theorem cr(CQ4)=8}, we conclude
that there is only one possible drawing of $\ell\mbox{-}LTQ_3$ as
shown in Figure 3.6(2).

Suppose $\cross(\langle\{a^{'},b^{'},c^{'},d^{'}\}\rangle)=0$  (see
Figure 3.5(3)). By \eqref{equation LTQ another two for four
4-cycles}, we may assume without loss of generality that
$$\{aa^{'},bb^{'},cd^{'},dc^{'}\}\subset E_{\ell}.$$ By
Claim \ref{claim two edges cross geq 1 or neq 1}, we have that
$$\cross(\{aa^{'},bb^{'}\}, E_{\ell}\setminus
\{aa^{'},bb^{'}\})+\cross(aa^{'},bb^{'})\geq 1,$$ and similarly,
$$\cross(\{cd^{'},dc^{'}\},E_{\ell}\setminus
\{cd^{'},dc^{'}\})+\cross(cd^{'},dc^{'})\neq 1.$$

If $\cross(\{cd^{'},dc^{'}\},E_{\ell}\setminus
\{cd^{'},dc^{'}\})+\cross(cd^{'},dc^{'})\geq 2$, then
$\cross(E_{\ell})\geq 3$. By Claim \ref{claim cr(El) geq 3}, we have
that there exists a region $f$ of $\ell\mbox{-}LTQ_3$ with
$|V_{on}(f)|\geq 7$. We may assume that $f$ is the unbounded region
of $\ell\mbox{-}LTQ_3$. It follows that there does not exist any
region $f$ of $\ell\mbox{-}LTQ_3$ with $|V_{on}(f)|=8$, otherwise,
$\cross(E_{\ell})=5$, a contradiction with \eqref{equation 1 leq
cr(El) leq 4}. By Claim \ref{claim cr(El) geq 3}, we infer that
$$\cross(E_{\ell})=3.$$ It follows from Claim \ref{claim two edges
cross geq 1 or neq 1} that $\cross(\{aa^{'},bb^{'}\}$,
$E_{\ell}\setminus \{aa^{'},bb^{'}\})+\cross(aa^{'},bb^{'})=1$ and
$\cross(\{cd^{'},dc^{'}\},E_{\ell}\setminus
\{cd^{'},dc^{'}\})+\cross(cd^{'},dc^{'})=2$. Combined with Claim
\ref{claim two cycles can not cross} and Claim \ref{claim four edges
are clean}, by symmetry, there is only one possible drawing of
$\ell\mbox{-}LTQ_3$ as shown in Figure 4.9.
\begin{figure}[h]
\centering
\includegraphics[scale=1.0]{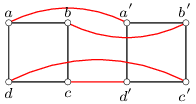}
\caption{\small{One possible drawing of $\ell\mbox{-}LTQ_3$}}
\end{figure}

If $\cross(\{cd^{'},dc^{'}\},E_{\ell}\setminus
\{cd^{'},dc^{'}\})+\cross(cd^{'},dc^{'})=0$, i.e., both $cd^{'}$ and
$dc^{'}$ are clean, by Claim \ref{claim four edges are clean}, we
have that all edges of $E(\{c,d,c^{'},d^{'}\})$ are clean. By Lemma
\ref{Lemma |Von(f)| leq 6}, we have that $|V_{on}(f)|\leq 6$ for
every region $f$ of $\ell\mbox{-}LTQ_3$. By \eqref{equation 1 leq
cr(El) leq 4} and Claim 1, we conclude that $\cross(E_{\ell})\leq
2$. Combined with Claim 2 and Claim 4, by symmetry, there is only
three possible drawings of $\ell\mbox{-}LTQ_3$ as shown in Figure
3.7(2)-(4).

Notice that the drawing shown in Figure 3.6(2) is isomorphic to the
drawing shown in Figure 3.7(3), and that the drawing shown in Figure
3.7(4) is isomorphic to the drawing shown in Figure 3.7(2). So we
need only to consider the drawings of $\ell\mbox{-}LTQ_3$ shown in
Figure 4.9 and Figure 3.7(2)-(3).

By Observation \ref{observation the paths P4 of LTQ_4}, we conclude
that there exist four vertex-disjoint 4-paths
$P_{au_a},P_{bu_b},P_{cu_c}$, $P_{du_d}$, where
$\{u_a,u_b\}=\{c^{'},d^{'}\}$ and $\{u_c,u_d\}=\{a^{'},b^{'}\}$. For
convenience, we denote $$E_p=E(P_{au_a})\cup E(P_{bu_b})\cup
E(P_{cu_c})\cup E(P_{du_d}).$$

\begin{claim} \label{claim the first drawing of l-LTQ3}
{\sl If $\ell\mbox{-}LTQ_3$ is drawn as Figure 4.9, then
$\cross(LTQ_4)\geq 10$.}
\end{claim}

{\sl Proof of Claim \ref{claim the first drawing of l-LTQ3}.} From
Figure 4.9, we see that
\begin{equation}\label{equation cr(El)=3}
\cross(E_{\ell})=3
\end{equation}
and
\begin{equation}\label{equation max(Von(f))=7}
|V_{on}(f)|\leq 7
\end{equation}
for every region $f$ of $\ell\mbox{-}LTQ_3$. By \eqref{left leq 3.5}
and \eqref{equation cr(El)=3}, we have
\begin{equation}\label{equation cr(El,Elr)+0.5cr(El,Er)leq 1.5(2)}
\cross(E_{\ell},E_{\ell,r})+0.5\cross(E_{\ell},E_r)\leq 1.5,
\end{equation}
which implies that $\cross(E_{\ell},E_r)\leq 3$. If
$\cross(E_{\ell},E_r)=3$, by \eqref{equation
cr(El,Elr)+0.5cr(El,Er)leq 1.5(2)} and Lemma \ref{Lemma the number
of emitting edges}, we conclude that all vertices of
$r\mbox{-}LTQ_3$ lie in exactly two regions $h_1, h_2$ of
$\ell\mbox{-}LTQ_3$ with $|V_{on}(h_1)|\geq |V_{on}(h_2)|$ and
$(|V_{in}(h_1;r\mbox{-}LTQ_3)|,|V_{in}(h_2;r\mbox{-}LTQ_3)|)=(7,1)$.
It follows that $h_1=f_1$ and $h_2\in\{f_2,f_4\}$. From Figure 4.7,
we see that the common boundary of $h_1$ and $h_2$ is an edge. By
Lemma \ref{lemma a good drawing about a vertex with degree 4}, this
is not a good drawing. Combined with \eqref{equation max(Von(f))=7}
and \eqref{equation cr(El,Elr)+0.5cr(El,Er)leq 1.5(2)}, we conclude
that
$$\cross(E_{\ell},E_r)=0 \mbox{ \ \ and \ \ } \cross(E_{\ell},E_{\ell,r})=1,$$
which implies that all vertices of $r\mbox{-}LTQ_3$ lie in the same
region $f_1$ of $\ell\mbox{-}LTQ_3$, and $b\pi(b)$ crosses $aa^{'}$
or $ad$.

Observe that every edge of $\{ab,dc,a^{'}b^{'},d^{'}c^{'}\}$ belongs
to exactly two 4-cycles of $\langle\{a,b,c,d\}\rangle$,
$\langle\{a^{'},b^{'},c^{'},d^{'}\}\rangle$,
$\langle\{a,b,a^{'},b^{'}\}\rangle$ and
$\langle\{c,d,c^{'},d^{'}\}\rangle$. By the structure of $LTQ_4$, we
have $v_0v_2\in\{ab,dc,a^{'}b^{'},d^{'}c^{'}\}$. Without loss of
generality, we may assume $a=v_0$ and $b=v_2$ (see Figure 4.10(1)
where $i$ stands for the vertex $v_i$ for $i\in\{0,1,\ldots,7\}$).
\begin{figure}[h]
\centering
\includegraphics[scale=1.0]{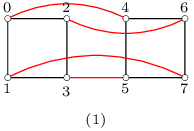} \hspace{30pt}
\includegraphics[scale=1.0]{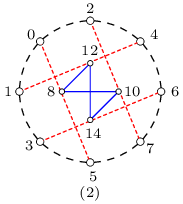} \hspace{30pt}
\includegraphics[scale=1.0]{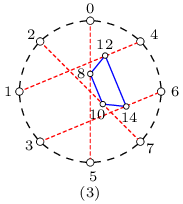}
\caption{\small{The graphs for the proof of Claim 7}}
\end{figure}

By the structure of $LTQ_4$ and Observation \ref{observation the
paths P4 of LTQ_4}, there exist four vertex-disjoint 4-paths
$P_{v_0v_5}=v_0v_8v_9v_5$, $P_{v_1v_4}=v_1v_{13}v_{12}v_{4}$,
$P_{v_2v_7}=v_2v_{10}v_{11}v_7$ and $P_{v_3v_6}=v_3v_{15}v_{14}v_6$.
Since $\cross(E_{\ell},E_{\ell,r})=1$, by \eqref{equation cr(LTQ4)
leq 9}, \eqref{equation cr(El)=3}, Lemma \ref{Lemma total crossing
and left } and Observation \ref{observation the relation of paths},
we conclude that $$\cross(E_p)\in\{4,5\},$$
$\cross(E(P_{v_0v_5}),E(P_{v_1v_4}))=1$,
$\cross(E(P_{v_0v_5}),E(P_{v_3v_6}))=1$,
$\cross(E(P_{v_2v_7}),E(P_{v_1v_4}))=1$,
$\cross(E(P_{v_2v_7}),E(P_{v_3v_6}))=1$,
$\cross(E(P_{v_1v_4}),E(P_{v_3v_6}))=0$ and
$\cross(E(P_{v_0v_5}),E(P_{v_2v_7}))\in\{0,1\}$.

If $\cross(E(P_{v_0v_5}),E(P_{v_2v_7}))=0$, then $\cross(E_p)=4$ and
$v_2v_{10}$ crosses $v_0v_4$ (see Figure 4.10(2)). Since $v_8\in
P_{v_0v_5}$, $v_{10}\in P_{v_2v_7}$, $v_{14}\in P_{v_3v_6}$ and
$v_{12}\in P_{v_1v_4}$, we see that at least one edge of the 4-cycle
$v_8v_{10}v_{14}v_{12}v_8$ is not clean and at least one edge of the
4-cycle $v_9v_{11}v_{13}v_{15}v_9$ is not clean. Hence,
$\cross(LTQ_4)\geq 10$, a contradiction with \eqref{equation
cr(LTQ4) leq 9}.

If $\cross(E(P_{v_0v_5}),E(P_{v_2v_7}))=1$, then $\cross(E_p)=5$ and
$v_2v_{10}$ crosses $v_0v_1$ (see Figure 4.10(3)). By
\eqref{equation cr(LTQ4) leq 9}, \eqref{equation cr(El)=3} and Lemma
\ref{Lemma total crossing and left }, we have $\cross(LTQ_4)=9$,
which implies that all edges of the two 4-cycles
$v_8v_{10}v_{14}v_{12}v_8$ and $v_9v_{11}v_{13}v_{15}v_9$ are clean.
From Figure 4.8(3), we see that $v_{11}v_{13}$ is not clean, a
contradiction. This proves Claim \ref{claim the first drawing of
l-LTQ3}. \qed

\begin{claim} \label{claim the second drawing of l-LTQ3}
{\sl If $\ell\mbox{-}LTQ_3$ is drawn as Figure 3.7(2), then
$\cross(LTQ_4)\geq 10$.}
\end{claim}

{\sl Proof of Claim \ref{claim the second drawing of l-LTQ3}.} For
convience, we lable all the regions of the drawing of Figure 3.7(2)
as shown in Figure 4.11.
\begin{figure}[h]
\centering
\includegraphics[scale=1.0]{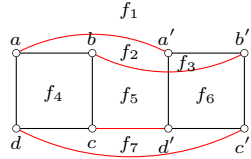}
\caption{\small{The graphs for the proof of Claim 8}}
\end{figure}
We see that
\begin{equation}\label{equation cr(El)=1}
\cross(E_{\ell})=1
\end{equation}
and
\begin{equation}\label{equation max(Von(f))=5}
|V_{on}(f)|\leq 5
\end{equation}
for every region $f$ of $\ell\mbox{-}LTQ_3$. By \eqref{left leq
3.5}, we have
\begin{equation}\label{equation cr(El,Elr)+0.5cr(El,Er)leq 3.5(2)}
\cross(E_{\ell},E_{\ell,r})+0.5\cross(E_{\ell},E_r)\leq 3.5.
\end{equation}
Let $h_1,\ldots,h_k$ be all the regions of $\ell\mbox{-}LTQ_3$ with
$|V_{in}(h_i;r\mbox{-}LTQ_3)|>0$ for $i=1,2,\ldots,k$. We shall
admit that
$$|V_{in}(h_1;r\mbox{-}LTQ_3)|\geq |V_{in}(h_2;r\mbox{-}LTQ_3)|\geq \cdots\geq
|V_{in}(h_k;r\mbox{-}LTQ_3)|.$$ By \eqref{equation
cr(El,Elr)+0.5cr(El,Er)leq 3.5(2)}, Lemma \ref{Lemma four region geq
7} and Lemma \ref{Lemma five regions geq 8}, we have
$$k\leq 4.$$

Suppose $k=4$. By \eqref{equation cr(El,Elr)+0.5cr(El,Er)leq 3.5(2)}
and Lemma \ref{Lemma four region geq 7}, we have that
$\cross(E_{\ell},E_r)=7$ and $\cross(E_{\ell},E_{\ell,r})=0$, which
implies
\begin{equation}\label{equation pi(Vin ) contained in Von}
\pi^{-1}(V_{in}(h_i;r\mbox{-}LTQ_3))\subseteq V_{on}(h_i)
\end{equation}
for $i=1,2,3,4$. By \eqref{equation pi(Vin ) contained in Von} and
Lemma \ref{Lemma four region geq 7}, we may assume without loss of
generality that $(h_1,h_2,h_3,h_4)\in
\{(f_1,f_2,f_4,f_7),(f_1,f_4,f_7,f_6)\}$. By \eqref{equation pi(Vin
) contained in Von} and Lemma \ref{Lemma four region geq 7}, we have
that
$\pi^{-1}(V_{in}(h_1;r\mbox{-}LTQ_3))=\{a,d,a^{'},b^{'},c^{'}\}$,
$\pi^{-1}(V_{in}(h_2;r\mbox{-}LTQ_3))=\{b\}$,
$\pi^{-1}(V_{in}(h_3;r\mbox{-}LTQ_3))=\{c\}$,
$\pi^{-1}(V_{in}(h_4;r\mbox{-}LTQ_3))=\{d^{'}\}$ and
$\pi(b)\pi(c),\pi(c)\pi(d^{'})\in E(LTQ_4)$. By Observation
\ref{observation the structure of P3 in l-LTQ3}, we derive a
contradiction. Therefore, $$k\leq 3.$$ \indent We consider the case
of $k=3$. For convenience, let $V_i=V_{in}(h_i;r\mbox{-}LTQ_3)$ for
$i=1,2,3$. By \eqref{equation cr(El,Elr)+0.5cr(El,Er)leq 3.5(2)} and
Lemma \ref{Lemma three vertex set (6 1 1)}, we conclude that $5\leq
\cross(E_{\ell},E_r)\leq 7.$

We first show that
\begin{equation}\label{equation cross of El and Er in [5,7]}
\cross(E_{\ell},E_r)\in\{6,7\}.
\end{equation}
Assume to the contrary that $\cross(E_{\ell},E_r)=5$. By Lemma
\ref{Lemma three vertex set (6 1 1)}, we have that
$$(|V_{in}(h_1;r\mbox{-}LTQ_3)|,|V_{in}(h_2;r\mbox{-}LTQ_3)|,
|V_{in}(h_3;r\mbox{-}LTQ_3)|)=(6,1,1)$$ and that $\mathcal
{B}(h_i,h_j)=1$ for any $(i,j)\in \{(1,2),(1,3),(2,3\}$. By
\eqref{equation max(Von(f))=5} and \eqref{equation
cr(El,Elr)+0.5cr(El,Er)leq 3.5(2)}, we have
$\cross(E_{\ell},E_{\ell,r})=1$ which implies $|V_{on}(h_1)\cup
V_{on}(h_2)\cup V_{on}(h_3)|\geq 7$. It follows that
$$\{h_1,h_2,h_3\}\in\{\{f_1,f_2,f_4\},\{f_1,f_4,f_7\},\{f_1,f_6,f_7\}\}.$$
Note that $P_{bu_b}$ crosses $aa^{'}$ or $ad$, and that $P_{cu_c}$
crosses  $ad$ or $dc^{'}$, and that $P_{au_a}$ crosses $dc^{'}$ or
$b^{'}c^{'}$ if $u_a=d^{'}$, or $P_{bu_b}$ crosses $dc^{'}$ or
$b^{'}c^{'}$ if $u_b=d^{'}$. By Observation \ref{observation the
relation of paths}, we conclude that $\cross(E_p)\geq 4$. It follows
that $\cross(LTQ_4)\geq
\cross(E_{\ell})+\cross(E_{\ell},E_r)+\cross(E_p)\geq 1+5+4=10$, a
contradiction with \eqref{equation cr(LTQ4) leq 9}. Therefore,
\eqref{equation cross of El and Er in [5,7]} holds.

By \eqref{equation cr(El,Elr)+0.5cr(El,Er)leq 3.5(2)} and
\eqref{equation cross of El and Er in [5,7]}, we have
\begin{equation}\label{equation cr(El,Elr)=0(2)}
\cross(E_{\ell},E_{\ell,r})=0,
\end{equation}
which implies  $$V_{on}(h_1)\cup V_{on}(h_2)\cup
V_{on}(h_3)=V(\ell\mbox{-}LTQ_3).$$ By \eqref{equation
max(Von(f))=5} and \eqref{equation cr(El,Elr)=0(2)}, we have that
\begin{equation}\label{equation Vin(fi) leq 5}
|V_{in}(h_i;r\mbox{-}LTQ_3)|\leq 5
\end{equation}
for $i=1,2,3$.

Next we show there exists  one pair $(s,t)$ of
$\{(1,2),(1,3),(2,3)\}$ such that
\begin{equation}\label{equation B(hi,hj)=1}
\mathcal {B}(h_i,h_j)=1
\end{equation}
 for any
$(i,j)\in\{(1,2),(1,3),(2,3)\}\setminus \{(s,t)\}$. Assume without
loss of generality to the contrary that $\mathcal
{B}(h_1,h_2)=\mathcal {B}(h_1,h_3)=0$. If
$|V_{in}(h_1;r\mbox{-}LTQ_3)|\geq 2$, by Lemma \ref{lemma counting
the number of crossings}, we have that $e(V_1,V_2)+e(V_1,V_3)\geq
4$, and thus, $\cross(E_{\ell},E_r)\geq 2(e(V_1,V_2)+e(V_1,V_3))=8$,
a contradiction with \eqref{equation cross of El and Er in [5,7]}.
If $|V_{in}(h_1;r\mbox{-}LTQ_3)|=1$, i.e.,
$e(V_1,V_2)+e(V_1,V_3)=3$, by Lemma \ref{Lemma three vertex set (6 1
1)}, we have that $e(V_1,V_2)+e(V_1,V_3)+e(V_2,V_3)\geq 5$, and
thus, $e(V_2,V_3)\geq 2$. It follows that $\cross(E_{\ell},E_r)=
2(e(V_1,V_2)+e(V_1,V_3))+e(V_2,V_3)\geq 8$, a contradiction with
\eqref{equation cross of El and Er in [5,7]}. Therefore,
\eqref{equation B(hi,hj)=1} holds.

Suppose that $\mathcal {B}(h_s,h_t)=0$. By \eqref{equation
cr(El,Elr)=0(2)}, we have that
$$\{h_1,h_2,h_3\}\in\{\{f_1,f_2,f_7\},\{f_1,f_4,f_6\}\}.$$
If $e(V_s,V_t)=0$, by Lemma \ref{lemma three vertex set geq 8}, we
have that there exists a region $h_i$ of $\ell\mbox{-}LTQ_3$ with
$|V_{in}(h_i;r\mbox{-}LTQ_3)|=1$, where $i=1,2,3$. Observe that the
common boundary of any two regions in Figure 4.11 is an edge. By
Lemma \ref{lemma a good drawing about a vertex with degree 4}, this
is not a good drawing. Hence
$$e(V_s,V_t)\geq 1.$$
If $e(V_s,V_t)\geq 2$, say $e(V_1,V_2)\geq 2$, since
$\cross(E_{\ell},E_r)=e(V_1,V_2)\times 2+e(V_1,V_3)+e(V_2,V_3)\leq
7$, we have $e(V_1,V_3)+e(V_2,V_3)\leq 3$. By Lemma \ref{Lemma the
number of emitting edges}, we have $e(V_1,V_3)+e(V_2,V_3)\geq 3$. It
follows that $e(V_1,V_3)+e(V_2,V_3)=3$, and thus $e(V_1,V_2)=2$. By
Lemma \ref{Lemma three vertex set (6 1 1)}, we conclude that there
exists a region $h_i$ of $\ell\mbox{-}LTQ_3$ with
$|V_{in}(h_i;r\mbox{-}LTQ_3)|=6$ where $i=1,2,3$, a contradiction
with \eqref{equation Vin(fi) leq 5}. Hence
$$e(V_s,V_t)=1.$$ If $\cross(E_{\ell},E_r)=6$, then
$e(V_1,V_2)+e(V_1,V_3)+e(V_2,V_3)=5$. By Lemma \ref{Lemma three
vertex set (6 1 1)}, we conclude that there exists a region $h_i$ of
$\ell\mbox{-}LTQ_3$ with $|V_{in}(h_i;r\mbox{-}LTQ_3)|=6$ where
$i=1,2,3$, a contradiction with \eqref{equation Vin(fi) leq 5}.
Hence $$\cross(E_{\ell},E_r)=7.$$ Then
$e(V_1,V_2)+e(V_1,V_3)+e(V_2,V_3)=6$. By Lemma \ref{lemma three
vertex set equal 6}, we have that
$$(|V_{in}(h_1;r\mbox{-}LTQ_3)|,|V_{in}(h_2;r\mbox{-}LTQ_3)|,|V_{in}(h_3;r\mbox{-}LTQ_3)|)\in\{(5,2,1),(4,3,1)\}.$$
Then $\cross(E(P_{du_d}),E(P_{au_a}))\geq 1$ and
$\cross(E(P_{cu_c}),E(P_{au_a}))\geq 1$, i.e., $\cross(E_p)\geq 2$.
It follows that $\cross(LTQ_4)\geq
\cross(E_{\ell})+\cross(E_{\ell},E_r)+\cross(E_p)\geq 1+7+2=10$, a
contradiction with \eqref{equation cr(LTQ4) leq 9}. Therefore,
$$\mathcal {B}(h_i,h_j)=1$$
for any $(i,j)\in\{(1,2),(1,3),(2,3)\}$.

By \eqref{equation cr(El,Elr)=0(2)}, we conclude that
$\{h_1,h_2,h_3\}=\{f_1,f_4,f_7\}$, $\pi(a^{'}),\pi(b^{'})\in
V_{in}(f_1;r\mbox{-}LTQ_3)$, $\pi(b)\in V_{in}(f_4;r\mbox{-}LTQ_3)$,
$\pi(d^{'})\in V_{in}(f_7;r\mbox{-}LTQ_3)$, and $\pi(c)\in
V_{in}(f_4;r\mbox{-}LTQ_3)$ or $\pi(c)\in
V_{in}(f_7;r\mbox{-}LTQ_3)$. By \eqref{equation max(Von(f))=5} and
\eqref{equation cr(El,Elr)=0(2)}, we conclude that
\begin{equation}\label{equation vin neq 6 1 1}
(|V_{in}(h_1;r\mbox{-}LTQ_3)|,|V_{in}(h_2;r\mbox{-}LTQ_3)|,|V_{in}(h_3;r\mbox{-}LTQ_3)|)\neq
(6,1,1).
\end{equation}
By Observation \ref{observation vertex partition of LTQ3}, we have
that $\langle\{\pi(a),\pi(b),\pi(c),\pi(d)\}\rangle\cong C_4$ and
$\langle\{\pi(a^{'}),\pi(b^{'}),\pi(c^{'})$,
$\pi(d^{'})\}\rangle\cong C_4$, or
$\langle\{\pi(a),\pi(b),\pi(a^{'}),\pi(b^{'})\}\rangle\cong C_4$ and
$\langle\{\pi(c),\pi(d),\pi(c^{'}),\pi(d^{'})\}\rangle\cong C_4$.
Let $$E^{'}=E_r\setminus E_p.$$ From Figure 4.11, we see that the
four vertices $a^{'},b^{'},c^{'},d^{'}$ do not lie in the same
region and the four vertices $a,b,a^{'},b^{'}$ do not lie in the
same region. It follows from \eqref{equation cr(El,Elr)=0(2)} that
\begin{equation}\label{equation cr(Er-Ep) geq 1}
\cross(E^{'})+\cross(E^{'},E(LTQ_4)\setminus E^{'})\geq 1.
\end{equation}
If $\cross(E_{\ell},E_r)=7$, by \eqref{equation vin neq 6 1 1} and
Lemma \ref{Lemma three vertex set (6 1 1)}, we conclude that there
exists exactly one edge of $E_r$ span a region other than
$h_1,h_2,h_3$ from $h_i$ into $h_j$, where $1\leq i<j\leq 3$. This
implies that $P_{cu_c}$ crosses at least one edge of
$\{ab,ad,dc,dc^{'},d^{'}c^{'}\}$. It follows that
$\cross(E(P_{au_a}),E(P_{cu_c}))\geq 1$ or
$\cross(E(P_{bu_b}),E(P_{cu_c}))\geq 1$, i.e., $\cross(E_p)\geq 1$.
Since $\cross(E_{\ell},E_r)=7$, it follows from \eqref{equation
cr(El)=1} and \eqref{equation cr(Er-Ep) geq 1} that
$\cross(LTQ_4)\geq
\cross(E_{\ell})+\cross(E_{\ell},E_r)+\cross(E_p)+\cross(E^{'})+\cross(E^{'},E(LTQ_4)\setminus
E^{'})=1+7+1+1=10$, a contradiction with \eqref{equation cr(LTQ4)
leq 9}. Hence, we have
\begin{equation}\label{equation cr(El,Er)=6}
\cross(E_{\ell},E_r)=6
\end{equation}
By \eqref{equation vin neq 6 1 1}, \eqref{equation cr(El,Er)=6} and
Lemma \ref{Lemma three vertex set (6 1 1)}, we conclude that the
edges in $E_r$ cross the common boundary of the three regions
$h_1,h_2$ and $h_3$ only once. This implies that $P_{cu_c}$ crosses
$ad$ or $dc^{'}$. It follows that
$\cross(E(P_{au_a}),E(P_{cu_c}))\geq 1$ and
$\cross(E(P_{bu_b}),E(P_{cu_c}))\geq 1$, i.e., $\cross(E_p)\geq 2$.
Combined with \eqref{equation cr(El)=1}, \eqref{equation cr(Er-Ep)
geq 1} and \eqref{equation cr(El,Er)=6}, we have that
$\cross(LTQ_4)\geq
\cross(E_{\ell})+\cross(E_{\ell},E_r)+\cross(E_p)+\cross(E^{'})+\cross(E^{'},E(LTQ_4)\setminus
E^{'})=1+6+2+1=10$, a contradiction with \eqref{equation cr(LTQ4)
leq 9}. This completes the arguments for the case of $k=3$.

We proceed to consider the case of $k=2$. We first show that
\begin{equation}\label{equation B(h1,h2)=1}
\mathcal {B}(h_1,h_2)=1.
\end{equation}
Assume to the contrary that $\mathcal {B}(h_1,h_2)=0.$ By
\eqref{equation cr(El,Elr)+0.5cr(El,Er)leq 3.5(2)} and Lemma
\ref{Lemma the number of emitting edges}, we have $6\leq
\cross(E_{\ell},E_r)<8$ and
$(|V_{in}(h_1;r\mbox{-}LTQ_3)|,|V_{in}(h_2;r\mbox{-}LTQ_3)|)=(7,1)$.
It follows from \eqref{equation max(Von(f))=5} that
$\cross(E_{\ell},E_{\ell,r})\geq 2$. Since $\cross(E_r)\geq 1$, it
follows that $\cross(LTQ_4)\geq
\cross(E_{\ell})+\cross(E_r)+\cross(E_{\ell},E_r)+\cross(E_{\ell},E_{\ell,r})=1+1+6+2=10$,
a contradiction with \eqref{equation cr(LTQ4) leq 9}. Hence,
\eqref{equation B(h1,h2)=1}  holds.

By Lemma \ref{Lemma the number of emitting edges}, we have
$\cross(E_{\ell},E_r)\geq 3$. If $\cross(E_{\ell},E_r)=3$, by
\eqref{equation max(Von(f))=5}, \eqref{equation
cr(El,Elr)+0.5cr(El,Er)leq 3.5(2)} and Lemma \ref{Lemma the number
of emitting edges}, we conclude that $\cross(E_{\ell},E_{\ell,r})=2$
and
$(|V_{in}(h_1;r\mbox{-}LTQ_3)|,|V_{in}(h_2;r\mbox{-}LTQ_3)|)=(7,1)$.
From Figure 4.11, we see that
$$\{h_1,h_2\}\in\{\{f_1,f_2\},\{f_1,f_4\},\{f_1,f_6\},\{f_1,f_7\}\}$$
and that the common boundary of $h_1$ and $h_2$ is an edge. By Lemma
\ref{lemma a good drawing about a vertex with degree 4}, this is not
a good drawing. Hence $$\cross(E_{\ell},E_r)\geq 4.$$ By
\eqref{equation cr(El,Elr)+0.5cr(El,Er)leq 3.5(2)}, we have
$\cross(E_{\ell},E_{\ell,r})\leq 1$, which implies that
$|V_{on}(h_1)\cup V_{on}(h_2)|\geq 7$. From Figure 5.9(2), we see
that $|V_{on}(h_1)\cup V_{on}(h_2)|\leq 7$. Combined with
\eqref{equation cr(El,Elr)+0.5cr(El,Er)leq 3.5(2)}, we have that
\begin{equation}\label{equation cr(El,Elr)=1}
\cross(E_{\ell},E_{\ell,r})=1
\end{equation}
and
\begin{equation}\label{equation cr(El,Er)=4 or 5}
\cross(E_{\ell},E_r)\in \{4,5\}.
\end{equation}
It follows that $$\{h_1,h_2\}\in\{\{f_1,f_4\},\{f_1,f_7\}\}.$$ From
Figure 4.11, we see that $\pi(\{a^{'},b^{'},c^{'}\})\subset
V_{in}(f_1;r\mbox{-}LTQ_3)$ and $\pi(\{b,c\})\subset
V_{in}(f_4;r\mbox{-}LTQ_3)$, or $\pi(\{a^{'},b^{'},a\})\subset
V_{in}(f_1;r\mbox{-}LTQ_3)$ and $\pi(\{c,d^{'}\})\subset
V_{in}(f_7;r\mbox{-}LTQ_3)$, which implies that
$|V_{in}(h_i;r\mbox{-}LTQ_3)|\geq 2$ for $i=1,2$. It follows that
$$(|V_{in}(h_1;r\mbox{-}LTQ_3)|,|V_{in}(h_2;r\mbox{-}LTQ_3)|)\in\{(6,2),(5,3),(4,4)\}.$$
\indent In the following, we consider only the case of
$\{h_1,h_2\}=\{f_1,f_4\}$. The case of $\{h_1,h_2\}=\{f_1,f_7\}$ can
be dealt in a similar argument.

Suppose $\pi(d^{'})\in V_{in}(f_1;r\mbox{-}LTQ_3)$. If
$\cross(E_{\ell},E_r)=4$, by Lemma \ref{Lemma the number of emitting
edges}, we infer that the edges in $E_r$ cross the common boundary
of $h_1,h_2$ only once. This implies that $P_{bu_b}$ crosses $ad$,
and that $P_{cu_c}$ crosses $ad$. It follows that $\cross(E_p)\geq
4$. It follows that $\cross(LTQ_4)\geq
\cross(E_{\ell})+\cross(E_{\ell},E_{\ell,r})+\cross(E_{\ell},E_r)+\cross(E_p)\geq
1+1+4+4=10$, a contradiction with \eqref{equation cr(LTQ4) leq 9}.
Hence, $$\cross(E_{\ell},E_r)=5.$$ By Lemma \ref{Lemma the number of
emitting edges}, we conclude that there exist only one edge of $E_r$
span one region $f_i$ from $h_1$ into $h_2$, where $1\leq i\leq 7$
and $i\not\in\{1,4,7\}$. This implies that one of the two paths
$P_{bu_b}$ and $P_{cu_c}$ crosses $ab$ and $aa^{'}$, or $dc$ and
$dc^{'}$. By Observation \ref{observation the relation of paths}, we
conclude that $\cross(E_p)\geq 3$. It follows that
$\cross(LTQ_4)\geq
\cross(E_{\ell})+\cross(E_{\ell},E_{\ell,r})+\cross(E_{\ell},E_r)+\cross(E_p)\geq
1+1+5+3=10$, a contradiction with \eqref{equation cr(LTQ4) leq 9}.

Suppose $\pi(d^{'})\in V_{in}(f_4;r\mbox{-}LTQ_3)$. Then
$$(|V_{in}(h_1;r\mbox{-}LTQ_3)|,|V_{in}(h_2;r\mbox{-}LTQ_3)|)\in
\{(4,4),(5,3)\}.$$ \indent If
$(|V_{in}(h_1;r\mbox{-}LTQ_3)|,|V_{in}(h_2;r\mbox{-}LTQ_3)|)=(4,4)$,
by \eqref{equation cr(El,Er)=4 or 5} and Lemma \ref{Lemma the number
of emitting edges}, we conclude that $\langle
V_{in}(f_4;r\mbox{-}LTQ_3)\rangle\cong C_4$. By Observation
\ref{observation vertex partition of LTQ4}(ii), we derive a
contradiction.

If
$(|V_{in}(h_1;r\mbox{-}LTQ_3)|,|V_{in}(h_2;r\mbox{-}LTQ_3)|)=(5,3)$,
by \eqref{equation cr(El,Er)=4 or 5} and Lemma \ref{Lemma the number
of emitting edges}, we conclude that $\cross(E_{\ell},E_r)=5$,
$\langle V_{in}(h_1;r\mbox{-}LTQ_3)\rangle\cong C_5$ and the edges
in $E_r$ cross the common boundary of $h_1,h_2$ only once. By Lemma
\ref{Lemma crossings above 4-cycle and 5-cycle}, we conclude that
$\cross(E_r)+\cross(E_r,E_{\ell,r})\geq 3$. It follows that
$\cross(LTQ_4)\geq
\cross(E_{\ell})+\cross(E_{\ell},E_{\ell,r})+\cross(E_{\ell},E_r)+\cross(E_r)+\cross(E_r,E_{\ell,r})\geq
1+1+5+3=10$, a contradiction with \eqref{equation cr(LTQ4) leq 9}.

Finally, we consider the case of $k=1$. By \eqref{equation
max(Von(f))=5} and \eqref{equation cr(El,Elr)+0.5cr(El,Er)leq
3.5(2)}, we conclude that $h_1=f_1$ and
$\cross(E_{\ell},E_{\ell,r})=3$, which implies that $P_{bu_b}$
crosses $aa^{'}$ or $ad$, and that $P_{cu_c}$ crosses $ad$ or
$dc^{'}$, and that $P_{au_a}$ crosses $dc^{'}$ or $b^{'}c^{'}$ if
$u_a=d^{'}$, or $P_{bu_b}$ crosses $dc^{'}$ or $b^{'}c^{'}$ if
$u_b=d^{'}$. Since $\cross(E_{\ell},E_{\ell,r})=3$, by
\eqref{equation cr(LTQ4) leq 9} and \eqref{equation cr(El)=1}, we
conclude that $ \cross(E_p)\in\{4,5\}$. By a similar argument as
Claim \ref{claim the first drawing of l-LTQ3}, we derive a
contradiction. This proves Claim \ref{claim the second drawing of
l-LTQ3}.\qed

\begin{claim} \label{claim the third drawing of l-LTQ3}
{\sl If $\ell\mbox{-}LTQ_3$ is drawn as Figure 3.7(3), then
$\cross(LTQ_4)\geq 10$.}
\end{claim}

{\sl Proof of Claim \ref{claim the third drawing of l-LTQ3}.} For
convience, we lable all the regions of the drawing of Figure 3.7(3)
as shown in Figure 4.12.
\begin{figure}[h]
\centering
\includegraphics[scale=1.0]{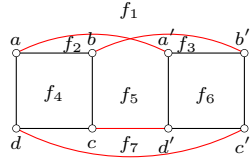}
\caption{\small{The graphs for the proof of Claim 9}}
\end{figure}

We see that $\cross(E_{\ell})=1$ and
\begin{equation}\label{equation max(Von(f))=4}
|V_{on}(f)|\leq 4
\end{equation}
for every region $f$ of $\ell\mbox{-}LTQ_3$. By \eqref{left leq
3.5}, we have
\begin{equation}\label{equation cr(El,Elr)+0.5cr(El,Er)leq 3.5}
\cross(E_{\ell},E_{\ell,r})+0.5\cross(E_{\ell},E_r)\leq 3.5.
\end{equation}
By \eqref{equation max(Von(f))=4}, \eqref{equation
cr(El,Elr)+0.5cr(El,Er)leq 3.5}, Lemma \ref{Lemma four region geq 7}
and Lemma \ref{Lemma five regions geq 8}, we conclude that
$$2\leq |\{h: h \mbox { is a region of }r\mbox{-}LTQ_3, V_{in}(h;r\mbox{-}LTQ_3)>0\}|\leq 3.$$
\indent Suppose that all vertices of $r\mbox{-}LTQ_3$ lie in exactly
three regions $h_1,h_2,h_3$ of $\ell\mbox{-}LTQ_3$. By
\eqref{equation max(Von(f))=4}, \eqref{equation
cr(El,Elr)+0.5cr(El,Er)leq 3.5} and Lemma \ref{Lemma three vertex
set (6 1 1)}, we conclude that $\cross(E_{\ell},E_r)\in\{6,7\}$ and
$\cross(E_{\ell},E_{\ell,r})=0$, which implies that $V_{on}(h_1)\cup
V_{on}(h_2)\cup V_{on}(h_3)=V(\ell\mbox{-}LTQ_3)$. From Figure
3.4(3), we see that there exists $(s,t)\in \{(1,2),(1,3),(2,3)\}$
such that $\mathcal {B}(h_s,h_t)=0$. Moreover, since
$\cross(E_{\ell},E_r)\in\{6,7\}$, we have that $$\mathcal
{B}(h_i,h_j)=1$$ for any $(i,j)\in \{(1,2),(1,3),(2,3)\}\setminus
\{(s,t)\}$. Then we observe that
$\{h_1,h_2,h_3\}\in\{\{f_1,f_2,f_5\}$,
$\{f_1,f_3,f_5\},\{f_1,f_4,f_5\},\{f_1,f_5,f_6\},\{f_1,f_5,f_7\},\{f_4,f_5,f_6\}\}$.
By a similar argument as Claim \ref{claim the second drawing of
l-LTQ3}, we derive a contradiction.

Suppose that all vertices of $r\mbox{-}LTQ_3$ lie in exactly two
regions $h_1,h_2$ of $\ell\mbox{-}LTQ_3$. By \eqref{equation
max(Von(f))=4} and Lemma \ref{Lemma the number of emitting edges},
we conclude that $\cross(E_{\ell},E_r)\geq 4$, $|V_{on}(h_1)\cup
V_{on}(h_2)|\geq 7$ and $\mathcal {B}(h_1,h_2)=1$. From Figure 4.12,
we see that there does not exist such two regions. This proves Claim
\ref{claim the third drawing of l-LTQ3}.\qed

Combined Claim \ref{claim the first drawing of l-LTQ3}, Claim
\ref{claim the second drawing of l-LTQ3} and Claim \ref{claim the
third drawing of l-LTQ3}, we derive a conradicition. This completes
the proof of Theorem \ref{theorem cr(LTQ4)=10}.
\end{proof}

\section{Conclusion}

\indent \indent In Section 3 and Section 4, we obtained that the
crossing numbers of $CQ_4$ and $LTQ_4$ are $8$ and $10$,
respectively. As for M\"{o}bius cube, since $0\mbox{-}MQ_4\cong
LTQ_4$, by Theorem \ref{theorem cr(LTQ4)=10}, we have
\begin{theorem}
$cr(0\mbox{-}MQ_4)=10$.
\end{theorem}

In Figure 5.1, we give a good drawing of $1\mbox{-}MQ_4$ with 10
crossings, hence the crossing number of  $1\mbox{-}MQ_4$ is no more
than 10.
\begin{figure}[h]
\centering
\includegraphics[scale=0.9]{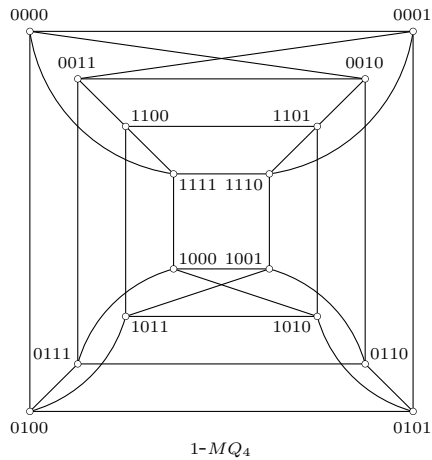}
\caption{\small{A good drawing of $1\mbox{-}MQ_4$ with 10
crossings}}
\end{figure}
Moreover, by an argument similar as the proof of $LTQ_4$ in Theorem
\ref{theorem cr(LTQ4)=10}, we have
\begin{theorem}
$cr(1\mbox{-}MQ_4)=10$.
\end{theorem}

As mentioned before, determination of the exact value of crossing
number for any kind of graph is a hard problem. For all kinds of
variations of hypercubes, all the determined value of crossing
number are summarized in Table 5.1, where the results obtained in
this paper are emphasized in bold fonts.

\begin{table}[htbp]
\centering
\begin{tabular}{|l|c|c|c|c|c}
  \hline
  $n$               & \ \ 1 \ \ & \ \ 2 \ \ & \ \ 3 \ \ & \ \ 4 \ \ \\ \hline
  $Q_n$             & 0 & 0 & 0 & 8 \\ \hline
  $CQ_n$            & 0 & 0 & 1 & $\textbf{8}$ \\ \hline
  $LTQ_n$           & 0 & 0 & 1 & $\textbf{10}$ \\ \hline
  $0\mbox{-}MQ_n$   & 0 & 0 & 1 & $\textbf{10}$ \\ \hline
  $1\mbox{-}MQ_n$   & 0 & 0 & 1 & $\textbf{10}$ \\ \hline
\end{tabular}
\caption{\small{The values of crossing numbers for the hypercube
variations}}
\end{table}

\end{document}